\documentclass[12pt]{amsart}
\usepackage{graphicx}
\usepackage{amsmath, enumerate}
\usepackage{amsfonts, amssymb, amsthm, xcolor, stmaryrd}
\usepackage{tikz, tikz-cd, caption, tabu, longtable, mathrsfs, xtab, centernot, mathtools}
\usepackage{dsfont, cite, todonotes}
\usepackage[toc,page]{appendix}
\usepackage[pagebackref=true]{hyperref}

\numberwithin{equation}{section}
\theoremstyle{plain}

\newtheorem{thm}{Theorem}[section]
\newtheorem{lemma}[thm]{Lemma}
\newtheorem{prop}[thm]{Proposition}
\newtheorem{cor}[thm]{Corollary}

\theoremstyle{definition}

\newtheorem{exmp}[thm]{Example}

\theoremstyle{remark}
\newtheorem{rmk}[thm]{Remark}

\newcommand{\Br}{}

\newcommand{\lcm}{\mathrm{lcm}}

\newcommand{\Gr}{\mathrm{Gr}}
\newcommand{\Rb}{\mathbb{R}}
\newcommand{\Zb}{\mathbb{Z}}
\newcommand{\Qb}{\mathbb{Q}}
\newcommand{\Cb}{\mathbb{C}}
\newcommand{\ef}{\mathfrak{e}}
\newcommand{\ebf}{{\mathbf{e}}}
\newcommand{\half}{{\frac{1}{2}}}
\newcommand{\reg}{{\mathrm{reg}}}
\newcommand{\Hb}{\mathbb{H}}
\newcommand{\Mp}{{\mathrm{Mp}}}
\newcommand{\SL}{{\mathrm{SL}}}

\newcommand{\Nb}{\mathbb{N}}
\newcommand{\smat}[4]{\left(\begin{smallmatrix}
                 #1 & #2\\
                 #3 & #4
\end{smallmatrix}\right)}

\newcommand{\lp}{\left (}
\newcommand{\rp}{\right )}

\newcommand{\Oc}{\mathcal{O}}

\newcommand{\Tr}{\mathrm{Tr}}
\newcommand{\sgn}{\mathrm{sgn}}

\newcommand{\Nm}{\mathrm{Nm}}
\renewcommand{\a}{\mathfrak{a}}
\newcommand{\af}{\mathfrak{a}}
\newcommand{\mf}{\mathfrak{m}}
\newcommand{\vf}{\mathfrak{v}}

\newcommand{\CT}{\mathrm{CT}}
\newcommand{\db}{\mathrm{db}}
\newcommand{\e}{\epsilon}
\newcommand{\Vc}{\mathcal{V}}

\begin{document}
\title[Mock modular forms]{Mock modular forms with integral Fourier coefficients}
\author{Yingkun Li and Markus Schwagenscheidt}
\address{Fachbereich Mathematik,
Technische Universit\"at Darmstadt, Schlossgartenstrasse 7, D--64289
Darmstadt, Germany}
\email{li@mathematik.tu-darmstadt.de}

\address{ETH Z\"urich Mathematics Dept., R\"amistrasse 101, CH-8092 Z\"urich, Switzerland}
\email{mschwagen@ethz.ch}
\date{\today}

\begin{abstract}
  In this note, we explicitly construct mock modular forms with integral Fourier coefficients by evaluating regularized Petersson inner products involving their shadows, which are unary theta functions of weights $\frac{1}{2}$ and $\frac{3}{2}$. 
  In addition, we also improve the known bounds for the denominators of the coefficients of mock modular forms whose shadows are holomorphic weight one cusp forms constructed by Hecke. 
\end{abstract}
\maketitle

\section{Introduction}

In his groundbreaking thesis \cite{ZwegersThesis}, Zwegers discovered that Ramanujan's mock theta functions could be completed to real-analytic modular forms by adding suitable non-holomorphic functions, whose images under the lowering operator are essentially complex conjugates of weight $\frac{3}{2}$ unary theta functions. This observation clarified the precise modularity properties of the mock theta functions and opened up new ways to study these functions. In the sense of Bruinier and Funke \cite{BF04}, these real-analytic completions are harmonic Maass forms of weight $k = \half$, whose images under the differential operator $\xi_{k} := 2iv^k \overline{\frac{\partial}{\partial\overline{\tau}}}$ are unary theta functions of weight $\frac{3}{2}$.
The holomorphic part and $\xi$-image of a harmonic Maass form are called a mock modular form and its shadow (see \cite{Za09}).

To establish the modular property of Ramanujan's classical mock theta functions, Zwegers gave three different constructions of these weight $\half$ harmonic Maass forms using Appell-Lerch sums, indefinite theta functions, and Fourier coefficients of meromorphic Jacobi forms, respectively. Bringmann and Ono \cite{BringmannOno} proved that generating series of certain partition numbers yield mock modular forms of weight $\frac{1}{2}$ whose shadows are linear combinations of unary theta functions.
In the case of weight $\frac{3}{2}$ mock modular forms with holomorphic theta functions as shadows, the prominent example is the generating series of Hurwitz class numbers that
appeared in the seminal work of Hirzebruch and Zagier \cite{HZ76}. Building on the work of Zwegers, the paper \cite{BFO09} also constructed such mock modular forms and related them to $q$-series.

For an arbitrary unary theta function of weight $\frac{1}{2}$ or $\frac{3}{2}$, Bruinier and the second author \cite{BS17} gave a different construction of its $\xi$-preimage using regularized theta lifts, and expressed the Fourier coefficients of the holomorphic part in terms of CM values of modular functions. From these constructions, one can apply the theory of complex multiplication to show that these Fourier coefficients are rational. A drawback of the construction from \cite{BS17} is the fact that it does not yield an explicit bound on the denominators of the coefficients.

In \cite{DMZ12}, Dabholkar, Murty and Zagier studied properties of mixed mock modular forms, which are products of mock modular forms and holomorphic modular forms. Using meromorphic Jacobi forms, they constructed mock modular forms whose shadows are unary theta functions and whose Fourier coefficients are rational numbers with explicitly bounded denominators (see equation (9.2) loc.\ cit.). However, the bound grows exponentially with respect to the index, even though it seems possible to find smaller bound when the index is small (see \cite[section 9.5.3]{DMZ12}). 

In the present article, we give another construction of $\xi$-preimages of weight $\frac{1}{2}$ and $\frac{3}{2}$ unary theta functions, and show that their Fourier coefficients are rational numbers with absolutely bounded denominators. Let us explain our results in more detail.

For $N \in \Nb$ we denote by $\Cb[\Zb/2N\Zb]$ the group ring spanned by the formal basis symbols $\ef_h$ for $h \in \Zb/2N\Zb$. Then for $\nu \in \{0,1\}$ the vector-valued unary theta function
\begin{align}\label{eq:thetaN}
\theta_{N}(\tau; \nu) := \sum_{h\bmod{2N}}\theta_{N, h}(\tau; \nu)\ef_h, \qquad \theta_{N, h}(\tau; \nu) := \sum_{n \in 2N\Zb + h}  {n}^\nu q^{ \frac{n^2}{4N}},
\end{align}
with $q:=e^{2\pi i \tau}$, is a holomorphic modular form of weight $\nu + \frac{1}{2}$ for the Weil representation $\rho_N$ of $\Mp_2(\Zb)$ (see Section~\ref{sec:preliminaries}). Our first main result is as follows.

\begin{thm}
  \label{thm:denom}
For $N, d \in \Nb$, denote $N_d := \gcd(N, d)$. 
  Let $\theta_{N}(\tau; \nu)$ be as in \eqref{eq:thetaN}. 
  \begin{enumerate}
  	\item There exists a mock modular form $\widetilde{\theta}^+_{N}(\tau;1)$ of weight $\frac{1}{2}$ with shadow $\frac{1}{\sqrt{N}} \theta_{N}(\tau; 1)$ such that $24N_4 \widetilde{\theta}^+_{N}(\tau;1)$ has integral Fourier coefficients.
	\item There exists a mock modular form $\widetilde{\theta}^+_{N}(\tau;0)$ of weight $\frac{3}{2}$ with shadow $\frac{\sqrt{N}}{2\pi} \theta_{N}(\tau; 0)$ such that $144N_4 \widetilde{\theta}^+_{N}(\tau;0)$ has integral Fourier coefficients.
 \end{enumerate}
  \end{thm}

\begin{rmk}
  \label{rmk:level}
  As usual, $\widetilde{\theta}^+_{N}(\tau;\nu)$ is the holomorphic part of a harmonic Maass form $\widetilde{\theta}_N(\tau;\nu)$ of weight $\frac{3}{2} - \nu$ for the complex conjugate $\overline{\rho}_N$ of the Weil representation $\rho_N$ on $\Mp_2(\Zb)$.
  In general, they will have exponential growth near the cusp. From our construction, the order of the principal part can be controlled by $q^{-1/6}$. 
\end{rmk}

\begin{rmk}
  \label{rmk:denom}
  Here is a trick to reduce the denominator. 
  Given two mock modular forms $f_1, f_2$ with the same shadow $g$ such that $M_j f_j$ has integral Fourier coefficients for $M_j \in \Nb$, we can find $a_j \in \Zb$ such that $a_1M_1 + a_2M_2 = \gcd(M_1,M_2)$. Then $\frac{a_1 M_1 f_1 + a_2 M_2 f_2}{\gcd(M_1, M_2)}$ is mock modular with shadow $g$, and the denominators of its coefficients are bounded by $\gcd(M_1,M_2)$.
\end{rmk}
\begin{rmk}
  \label{rmk:optimal}
  The bound in Theorem \ref{thm:denom} could be reduced for certain $N$ (see e.g.\ Remark \ref{rmk:3N}). For small values of $N$, it certainly can be lowered by considering well-known examples. In general, perhaps a factor of 3 and small power of 2 can still be removed from case (2).
\end{rmk}
By adding up the components of the vector-valued mock modular form and scaling the variable $\tau$, we have the following result for scalar-valued unary theta functions as a direct consequence of Theorem \ref{thm:denom}.
\begin{cor}
  \label{cor:scalar}
  For  $\kappa \in \{0, 1\}, N \in \Nb$ and any periodic function $\varphi: \Zb/N\Zb\to \Zb$, there exists a mock modular form $\tilde\theta^+(\tau; \varphi, \kappa)$ of weight $3/2 - \kappa$ and level $4N$ such that $24\cdot 6^{1-\kappa} N_4 \tilde\theta^+(\tau; \varphi, \kappa)$ has integral Fourier coefficients, and the shadow of $ \tilde\theta^+(\tau; \varphi, \kappa)$ is the scalar-valued unary theta function
  $$
\theta(\tau; \varphi, \kappa) :=   (4\pi)^{\kappa-1} \sum_{n \in \Zb}\varphi(n)n^\kappa q^{n^2}
  $$
  of weight $\kappa + 1/2$ and level $4N$.
\end{cor}

The denominator bounds we obtained above can be used to give a denominator bound of the Weyl vectors for Borcherds products in signature $(1, 2)$ (see \cite{BS17}), which is closely related to the orders of the multiplier systems of such Borcherds products.\footnote{The order of the  multiplier system of a Borcherds product is always finite (see \cite{Bo00}).}
Also, the rational Fourier coefficients in Theorem \ref{thm:denom} appear in other formulas, such as the Fourier coefficients of mock modular forms with binary theta functions as shadows \cite{Ehlen17}, and CM values of higher Green's functions \cite{Li18, BEY19}. 
A denominator bound as we have proved here would make precise the fields of definition of the algebraic values appearing in those formulas.

For the proof of Theorem~\ref{thm:denom} we will compute, for any weakly holomorphic modular form $g$ of weight $\nu+\frac{1}{2}$ for $\rho_N$, the regularized Petersson inner product 
\[
(g(\tau),\theta_N(\tau;\nu))^\reg := \lim_{T \to \infty}\int_{\mathcal{F}_T} \langle g(\tau),\theta_N(\tau;\nu) \rangle v^{\nu + \frac{1}{2}} \frac{du dv}{v^2}, \qquad (\tau = u+iv \in \Hb),
\]
where $\mathcal{F}_T$ is the usual fundamental domain for $\SL_2(\Zb)$ truncated at the height $T$, and $\langle \cdot,\cdot \rangle$ denotes the natural hermitian pairing on $\Cb[\Zb/2N\Zb]$ (see Section~\ref{sec:preliminaries}). From this we will reconstruct $\widetilde{\theta}_{N}^+(\tau;\nu)$ using Serre duality (see Proposition \ref{prop:xipreimage}). The basic idea to compute the regularized Petersson inner product is to use the identity
\[
\langle g(\tau),\theta_{N}(\tau;\nu) \rangle =  \langle g(\tau)\eta(\tau)^{-4},\theta_{N}(\tau;\nu) \overline{\eta(\tau)^4}\rangle
\]
and realize $\theta_{N}(\tau;\nu) \overline{\eta(\tau)^4}$ as a special value of a signature $(1, 4)$ theta function. Then the inner product $(g(\tau),\theta_N(\tau;\nu))^\reg$ can be viewed as a special value of a regularized theta lift of the weakly holomorphic modular form $g(\tau)\eta(\tau)^{-4}$, which can be evaluated using the methods developed by Borcherds~\cite{Borcherds98} and Bruinier~\cite{bruinierhabil}. To make the argument rigorous, we will implement this idea for vector-valued forms.
This technicality causes us to have the factors $24N_4$ and $144N_4$ in the theorem above, which could  be reduced in certain cases (see Remark~\ref{rmk:3N}). We refer the reader to Section~\ref{sec:proof theorem denom} for the details of the proof of Theorem~\ref{thm:denom}.
Note that this idea has been used when $\theta_N(\tau; \nu)$ is replaced by a holomorphic binary theta function to produce harmonic Maass forms of weight one \cite{Ehlen17, Viazovska19}.

 The mock modular form $\widetilde{\theta}_{N}(\tau;\nu)$ from Theorem~\ref{thm:denom} will be constructed explicitly in the proof of the theorem, though it is not canonical and involves choosing a suitable lattice $L$, a primitive isotropic vector  $\ell \in L$ and its dual $\ell' \in L'$ (see e.g.\ \eqref{eq:ttheta}). By slightly modifying the method above, we will also give simpler formulas for mock modular forms with shadow $\theta_{N}(\tau; \nu)$ in Section~\ref{sec:explicit construction}. To give an impression, we include a special case in the introduction.
 \begin{prop}
   \label{prop:Hgen}
   Let $H(n)$  be the  Hurwitz class numbers   (with $H(0) := -\frac{1}{12}$). Then 
        \begin{equation}
          \label{eq:Hgen}
          \sum_{n \geq 0}H(n)q^{{n}{}}
          = \frac{1}{24\eta(4\tau)}\sum_{\mathfrak{a} \subset \Zb[\sqrt{6}]} \phi_6(\af) q^{\Nm(\af)/6}
          = \frac{1}{24\eta(4\tau)^3}\sum_{\mathfrak{a} \subset \Zb[\sqrt{2}]} \phi_2(\af) q^{\Nm(\af)/2}           
        \end{equation}
        with the functions $\phi_6, \phi_2$ defined by
        \begin{equation}
          \label{eq:phis}
          \begin{split}
            \phi_6(\af) &:=
            \begin{cases}
              \left( \frac{12}{\Tr_{}(\lambda)/2}\right)\Tr_{}\left(\frac{|\lambda|}{2-\sqrt{6}}\right)&    \text{ if }  \af = (\lambda) \text{ satisfies }
              49 - 20\sqrt{6} < \frac{\lambda}{\lambda'} \leq 1,\\
              0& \text{ otherwise,}
            \end{cases}\\
            \phi_2(\af) &:=
            \begin{cases}
              - \left( \frac{-4}{\Tr_{}(|\lambda|/2)}\right)\Tr_{}\left(\frac{\lambda^2\sqrt{2}}{3-\sqrt{2}} \right) & 
            \text{ if }  \af = (\lambda) \text{  satisfies }
            17-12\sqrt{2} < \frac{\lambda}{\lambda'} \leq 1,\\
            0& \text{ otherwise.}
            \end{cases}
          \end{split}
        \end{equation}
 \end{prop}
 
 The proposition follows from Propositions
 \ref{prop:mock theta weight 3half} and
\ref{prop:mock theta weight 3half alternative}
 below. There we construct for any $N \in \Nb$ an explicit mock modular form with shadow $\frac{\sqrt{N}}{\pi}\theta_N(\tau;0)$. Specializing to $N = 1$, and using that the generating series of Hurwitz class numbers is a mock modular form with shadow $-\frac{1}{8\pi}\theta_1(\tau;0)$ (see \cite{zagiereisenstein}), we easily obtain Proposition~\ref{prop:Hgen} by comparing the first few Fourier coefficients. The first identity can probably be derived from the classical results in \textsection 2.3 of \cite{HZ76}, whereas the second seems new to us.
Note that for $N = 2$, Proposition \ref{prop:mockthetaweightonehalf} also produces the following well-known mock modular form (see \cite{DMZ12})
\begin{equation}
  \label{eq:mocketa3}
    \begin{split}
  \widetilde{\theta}^+_2(\tau;1) &:=     \frac{E_2(\tau)/24 - F^{(2)}_2(\tau)}{\eta(\tau)^3} =
  -\frac{q^{-1/8}}{24} ( -1 + 45 q + 231q^2 + 770q^3 + O(q^4)),\\
  E_2(\tau) &:= 1 - 24\sum_{n \ge 1} \sigma_1(n) q^n, \qquad
  F^{(2)}_2(\tau) := \sum_{\substack{b > a > 0 \\ b - a \text{ odd}}} a (-1)^b q^{ab/2},
  \end{split}
\end{equation}
whose shadow is {$\frac{1}{2\sqrt{2}} \eta(\tau)^3$. }

  Using the same idea as in the proof of Theorem~\ref{thm:denom}, we can also improve a denominator bound for certain weight one harmonic Maass forms.
  Let $\mf \subset \Oc_F$ be an integral ideal in a real quadratic field $F = \Qb(\sqrt{D})$ with ring of integers $\Oc_F$, and   $\varphi$ an odd ray class group character of conductor $\mf \cdot \infty_1$ with $\Nm(\mf) = M$, which we view as a function on integral ideals by extending with 0.
  Hecke associated the holomorphic weight one, level $N = DM$ eigenform
  $$
f_\varphi(\tau) := \sum_{\af \subset \Oc_F} \varphi(\af) q^{\Nm(\af)}
  $$
  to $\varphi$ in \cite{Hecke26}.
  In \cite{CL20}, Charollois and the first author showed that there is a mock modular form $f_\varphi^+(\tau) = \sum_{n \gg -\infty} c_\varphi^+(n) q^n$ with shadow $f_\varphi$ satisfying
  \begin{equation}
  \label{eq:c+}
  c^+_\varphi(n) - 
\sum_{[\af] \in \mathrm{Cl}_F} \varphi(\af) \sum_{\substack{(\lambda) \subset \af \\ \Nm((\lambda) \af^{-1}) = n }} \overline{\varphi}(\lambda)  \log \left| \frac{\lambda}{\lambda'} \right|
\in 
\frac{1}{\kappa} \Zb[\varphi] \cdot \log \varepsilon_F
\end{equation}
with $\kappa = 96 D M^3 \prod_{p \mid 2M \text{ prime}} (1 + p^{-1})$.\footnote{There is a mistake in (4.2.7) of \cite{CL20}, where $A^3$ should be replaced by $A^2 D$.
This affects Theorem 1.1, 5.1, 6.5, where $M^3\phi(2M)$ should be replaced by $D M^2\phi(2M)$. }
Here $\varepsilon_F \in \Oc_F^\times$ is the fundamental unit, $\Zb[\varphi]\subset\Cb$ is the subring generated by the values of $\varphi$, and we have chosen a set of representatives of the class group $\mathrm{Cl}_F$.
The constant $\kappa$ comes from bounding the denominator of a mixed mock modular form of weight 1 in Theorem 5.1 of \cite{CL20}.
We can now improve it as follows.
\begin{thm}
  \label{thm:wt1}
In the notations above, we can take $\kappa = 96 N_4 = 96 (DM)_4 \mid 384$. 
  \end{thm}

  \begin{exmp}
     \label{exmp:eta2}
    The constant $96N_4$ above can sometimes be reduced. 
    Let $D= 12, \mf  = 2\sqrt{3} \Oc_F$ and define $\varphi(\af) = \sgn(\lambda) \cdot \varphi_0(\lambda)$ for an integral ideal $\af = (\lambda)$ with  $\varphi_0: \Oc_F \to \{ \pm 1\}$ given by
\begin{equation}
  \label{eq:psi0}
    \varphi_0(\lambda) := 
\begin{cases}
1 & \text{ if }\lambda \equiv 1, 2 + \sqrt{3} \pmod{2\sqrt{3}}, \\
-1 & \text{ if }\lambda \equiv 5, 4 + \sqrt{3} \pmod{2\sqrt{3}}, \\
0 & \text{ otherwise.}
\end{cases}
\end{equation}
Then $f_\varphi(\tau/12) = \sum_{\af \subset \Oc_F} \varphi(\af) q^{\Nm(\af)/12} =\eta(\tau)^2$.
Furthermore, we define
\begin{equation}
  \label{eq:tvarphi}
  \widetilde{\varphi}(\af) := \varphi(\af)  \log \left| \frac{\lambda}{\lambda'}\right|
\end{equation}
if
$ \af = (\lambda)$  with $2 - \sqrt{3} < \left| \frac{\lambda}{\lambda'} \right| \le 2 + \sqrt{3}$.
Then John Duncan asked if the function
\begin{equation}
  \label{eq:tvartheta}
  \widetilde{\vartheta}^+(\tau) := 4\log (2 + \sqrt{3})
  \eta(\tau) \widetilde{\theta}^+_2(\tau;1) 
  + \sum_{\af \subset \Oc_F} \widetilde{\varphi}(\af) q^{\Nm(\af)/12}
\end{equation}
is a mock modular form with shadow $\eta(\tau)^2$.
This is indeed the case and follows from the results in \cite{CL20}.
We will give the details in Section \ref{sec:proof theorem weight one}.
\end{exmp}

The paper is organized as follows. We start with a section on the necessary preliminaries about vector-valued harmonic Maass forms for the Weil representation, unary theta functions and their connection to the Dedekind eta function, and the relation between the evaluation of regularized inner products and mock modularity. In Section~\ref{sec:theta lifts} we evaluate several regularized theta lifts in signature $(1,n)$. Sections \ref{sec:proof theorem denom} and \ref{sec:proof theorem weight one} are devoted to the proofs of Theorem~\ref{thm:denom} and Theorem~\ref{thm:wt1}. Finally, in Section~\ref{sec:explicit construction} we give explicit constructions of mock modular forms whose shadows are unary theta functions.

\noindent \textbf{Acknowledgement}: We thank John Duncan for sharing the observation of Example \ref{exmp:eta2}.
We also thank an anynomous referee for bringing our attention to section 9 of \cite{DMZ12}, which led to improvement of the bound in Theorem \ref{thm:denom} from an earlier version.

The first author was supported by LOEWE research unit USAG, and by the Deutsche Forschungsgemeinschaft (DFG) through the Collaborative Research Centre TRR 326 "Geometry and Arithmetic of Uniformized Structures", project number 444845124.
The second author was supported by SNF project 200021\_185014.

\section{Preliminaries}
\label{sec:preliminaries}

\subsection{Modular forms for the Weil representation} We recall some facts about harmonic weak Maass forms for the Weil representation associated with an even lattice from \cite{bruinierhabil,BF04}.
Let ${\Mp}_2(\Rb)$ be the metaplectic two-fold cover of $\SL_2(\Rb)$ consisting of elements $(A,  \phi)$ with $ A = \smat{a}{b}{c}{d} \in \SL_2(\Rb)$ and $\phi: \Hb \to \Cb$ holomorphic with $\phi(\tau)^2 = c\tau + d$. Let $\Mp_2(\Zb)  \subset \Mp_2(\Rb)$ denote the inverse image of $\SL_2(\Zb)$ under the covering map $\Mp_2(\Rb) \to \SL_2(\Rb)$. It is generated by $T := (\smat{1}{1}{0}{1}, 1)$ and $S := (\smat{0}{-1}{1}{0}, \sqrt{\tau})$. We let $\widetilde{\Gamma}_\infty$ be the subgroup generated by $T$.

Let $L$ be an even lattice with quadratic form $Q$ of signature $(b^+, b^-)$, the dual lattice $L'$ and the associated finite quadratic module (or discriminant form)
\begin{equation}
  \label{eq:AL}
 A_L := L'/L,
\end{equation}
on which $Q$ becomes a quadratic form valued in $\Qb/\Zb$. Let $\Br(\cdot, \cdot)$ be the associated bilinear form. Moreover, we let $\Gamma_L$ be the discriminant kernel of $L$, which is the subgroup of the orthogonal group $\mathrm{O}(L)$ which acts trivially on $A_L$.
Let $\Cb[A_L] := \bigoplus_{h \in A_L} \Cb \ef_h$ be the group ring generated by $A_L$. 
Given a finite index sublattice $M \subset L$, we have a natural map $\varpi: M'/L \to A_L$, which induces a map
\begin{equation}
  \label{eq:sublat}
  \Cb[A_L] \to \Cb[A_M],~ \ef_h \mapsto \sum_{\delta \in A_M,~ \varpi(\delta) = h} \ef_\delta.
\end{equation}
For convenience, we use $L^-$ to denote the lattice $L$ with the quadratic form $-Q$. 

The group ring $\Cb[A_L]$ is naturally an $\Mp_2(\Zb)$-module via the Weil representation $\rho_L$ defined by
\begin{equation}
  \label{eq:Weil_rep}
  \begin{split}
    \rho_L(T)(\ef_h) &:= \ebf(Q(h)) \ef_h, \\
\rho_L(S)(\ef_h) &:= \frac{\ebf((b^- - b^+)/8)}{\sqrt{|A_L|}} \sum_{\mu \in A_L} \ebf(-\Br(h, \mu)) \ef_\mu,
  \end{split}
\end{equation}
where we put $\ebf(x):=e^{2\pi i x}$ for $x \in \Cb$. Despite the notation, $\rho_L$ only depends on the finite quadratic module $A_L$. 
There is a natural hermitian pairing $\langle \cdot , \cdot \rangle$ on $\Cb[A_L]$ given by $\langle \ef_{h_1}, \ef_{h_2} \rangle := 1$ if $h_1 = h_2$ and zero otherwise. With respect to this pairing $\rho_L$ is a unitary representation.

A real-analytic function $f: \Hb \to \Cb[A_L]$ is called modular of weight $k \in \half \Zb$ with respect to $\rho_L$ if 
\begin{equation}
  \label{eq:modular}
(f |_{k,L} (A, \phi))(\tau) := \phi(\tau)^{-2k}\rho_L^{-1}((A,\phi)) f(A z) = f(\tau)
\end{equation}
for all $(A, \phi) \in \Mp_2(\Zb)$ and $\tau \in \Hb$.
We denote the spaces of harmonic Maass, weakly holomorphic, holomorphic, and cuspidal forms of weight $k$ for $\rho_L$ by
\[
H_{k, L} \supset M^!_{k,L} \supset  M_{k, L} \supset S_{k, L},
\]
respectively. More generally, for any representation $\rho$ of $\Mp_2(\Zb)$, the subscript $L$ in the notation above will be replaced by $\rho$.

Every $f \in H_{k,L}$ can be written uniquely as $f = f^+ + f^-$, where $f^+$ is holomorphic and has a Fourier expansion of the form
\begin{align}\label{eq holomorphic part}
f^+(\tau) = \sum_{h \in A_L}\sum_{\substack{m \in \Zb+Q(h) \\ m \gg -\infty}}a_f(h,m)q^m \ef_h
\end{align}
with coefficients $a_f(h,m) \in \Cb$. We may assume that $2k \equiv \sgn(L) \pmod 2$ since otherwise the action $\rho_L(S^2)\ef_h = i^{-\sgn(L)}\ef_{-h}$ implies that $H_{k,L}$ is trivial. Moreover, under this assumption the coefficients above satisfy the symmetry 
\[
a_f(-h,m) = (-1)^{k-\frac{\sgn(L)}{2}}a_f(h,m)
\]
for every $h \in A_L$ and $m \in \Zb + Q(h)$. 

 The antilinear differential operator $\xi_{k} = 2iv^k \overline{\frac{\partial}{\partial \overline{\tau}}}$ defines a surjective map $H_{k,L} \to M_{2-k,L^-}^!$ (see Theorem~3.7 in \cite{BF04}). A holomorphic $\Cb[A_L]$-valued $q$-series as in \eqref{eq holomorphic part} is called a mock modular form of weight $k$ with shadow $g$ if it is the holomorphic part $f^+$ of a harmonic Maass form $f \in H_{k,L}$ which satisfies $g =\xi_k f$. 

Examples of harmonic Maass forms can be constructed as special values of Maass Poincar\'e series. Let $k \leq 0$. Following Section~1.3 in \cite{bruinierhabil}, for $h \in A_L$ and $m \in \Zb+Q(h)$ with $m < 0$ we consider the Maass Poincar\'e series
\[
F_{h,m}(\tau,s) := \frac{1}{2\Gamma(2s)}\sum_{(A,\phi) \in \widetilde{\Gamma}_\infty \backslash \Mp_2(\Zb)}\mathcal{M}_s(4\pi |m| v)\ebf(mu)\ef_{h}|_{k,L}(A,\phi),
\]
where $\mathcal{M}_s(v) := v^{-\frac{k}{2}}M_{-\frac{k}{2},s-\frac{1}{2}}(v)$ with the usual $M$-Whittaker function. We put\footnote{Note that for $k = 0$ the Maass Poincar\'e series $F_{h,m}(\tau,s)$ does not converge at $s = 1$, but it can be analytically continued to $s = 1$ via its Fourier expansion.}
\begin{align}\label{eq:maass poincare series}
F_{h,m}(\tau) := F_{h,m}\left(\tau,1-\frac{k}{2}\right) = q^{m}(\ef_h + \ef_{-h}) + O(1),
\end{align}
which defines a harmonic Maass form in $H_{k,L}$ that maps to a cusp form under $\xi_k$. 

\subsection{Unary theta series and the eta function}
\label{subsec:unary}

For $N \in \Nb$ we consider the lattice 
\begin{equation}
  \label{eq:ZN}
\Zb[N]:=(\Zb,Nx^2). 
\end{equation}
We can identify its discriminant form with $\Zb/2N\Zb$, and by a slight abuse of notation, we put $A_N := \Cb[\Zb/2N\Zb]$ and write $\rho_N$ for the Weil representation associated with $\Zb[N]$. The unary theta function $\theta_N(\tau;\nu)$ defined in \eqref{eq:thetaN} is a holomorphic modular form of weight $\nu+\frac{1}{2}$ for $\rho_N$. For $\nu = 1$ it is a cusp form.

Familiar modular forms can be expressed in terms of unary theta functions. For example,
\begin{equation}
  \label{eq:etatheta}
  \begin{split}
\eta(\tau) &=\theta_{6, 1}(\tau; 0) 
-\theta_{6, 5}(\tau; 0) ,\qquad 
 \eta(\tau)^3 = 
   \theta_{2, 1}(\tau;1).
  \end{split}
\end{equation}
To generalize this situation, it is convenient to phrase these identities in terms of eigenvectors of the Weil representation.
Define
\begin{equation}
  \label{eq:vm}
  \begin{split}
  \vf_3 &:=
  \ef_{(0, 1)} + \ef_{(1, 0)}
  + \ef_{(0, -1)} + \ef_{(-1, 0)}
  - (\ef_{(3, 2)} + \ef_{(2, 3)}
  + \ef_{(3, -2)} + \ef_{(-2, 3)})
  \in \Cb[A_3 \oplus A_3],\\
      \vf_2&:= \ef_1 - \ef_3 \in \Cb[A_2], \quad
  \vf_6:= \ef_1 - \ef_5 - \ef_7 + \ef_{11} \in \Cb[A_6], \quad
  \vf_4 := \vf_3 \otimes \vf_3 \in \Cb[A_3^4].
  \end{split}
\end{equation}
Here we write $(\sum_{\mu}a_\mu \ef_\mu) \otimes (\sum_{\nu}b_\nu \ef_\nu) := \sum_{\mu,\nu}a_\mu b_\nu \ef_{(\mu,\nu)}$.

\begin{lemma}
  \label{lemma:eigenvector}
  In the notations above, the vector $\vf_2$ (resp.\ $\vf_3, \vf_6, \vf_4$) generates a 1-dimensional subspace invariant under $\rho_2$ (resp.\ $\rho_3 \otimes \rho_3$, $\rho_6, \rho_3^{\otimes 4}$), which acts on the space via $\chi^3$ (resp.\ $\chi^2, \chi, \chi^4$). Here $\chi:  \Mp_2(\Zb) \to \Cb^\times$ is the character of $\Mp_2(\Zb)$ defined by
  \begin{equation}
  \label{eq:chi}
   	\chi(T) =  e^{2\pi i / 24}, \quad
    \chi(S) = 1.
  \end{equation}
  Furthermore, we have
  \begin{equation}
    \label{eq:etaid}
    \begin{split}
      \langle \theta_6(\tau; 0), \vf_6 \rangle &= 2 \eta(\tau), \\
      \langle \theta_2(\tau; 1), \vf_2 \rangle &= 2 \eta(\tau)^3, \\
      \langle \theta_3(\tau;0)^2      , \vf_3 \rangle &= 4 \eta(\tau)^2, \\
      \langle \theta_3(\tau;0)^4      , \vf_4 \rangle &= 16 \eta(\tau)^4.      
    \end{split}
  \end{equation}
\end{lemma}
\begin{rmk}
  Note that the first two equations in \eqref{eq:etaid} are equivalent to \eqref{eq:etatheta}, and the third equation has been known to Weber.
\end{rmk}

\begin{rmk}
  \label{rmk:chir}
  For all $0 \le r \le 23$, the space $M_{r/2, \chi^r} = S_{r/2, \chi^r}$ is 1-dimensional and spanned by $\eta(\tau)^r$. 
\end{rmk}
\begin{proof}[Proof of Lemma~\ref{lemma:eigenvector}]
  The first claim can be checked locally at each prime $p$, as
  \begin{align*}
      \vf_6 &=  (\ef_1 - \ef_{-1})  \otimes  (\ef_1 - \ef_{-1} )
  \in \Cb[A_6 \otimes \Zb_2] \otimes  \Cb[A_6 \otimes \Zb_3],\\
  \vf_3 &= (\ef_{(0, 1)} - \ef_{(1, 0)})
  \otimes
        (
        \ef_{(0, 1)}
        - \ef_{(1, 0)}
        + \ef_{(0, -1)}
        - \ef_{(-1, 0)}
        )  
\in   \Cb[A_3^2 \otimes \Zb_2] \otimes    \Cb[A_3^2 \otimes \Zb_3].
  \end{align*}
  The second claim follows from Remark \ref{rmk:chir} and comparing Fourier coefficients.
\end{proof}

\subsection{Regularized inner product and pairing.}
\label{subsec:pair}
  Let $L$ be an even lattice and $k \in \half \Zb$ satisfying $2k \equiv \sgn(L) \pmod{2}$. Denote
  \begin{equation}
    \label{eq:VA}
    V_{k,L} := \left\{
    \sum_{h \in A_L}\sum_{\substack{m \in \Zb+Q(h) \\ m \gg -\infty }} a(h, m) q^m \ef_h:
       a(-h, m) = (-1)^{k-\frac{\sgn(L)}{2}} a(h, m) \text{ for all } h, m
    \right\},
  \end{equation}
a vector space of formal Laurent series with values in $\Cb[A_L]$, which contains the space $M^!_{k, L}$ of weakly holomorphic modular forms. Furthermore, let $V^+_{k,L}$ denote the subspace consisting of formal power series (supported on indices $m \geq 0$), which contains the space $M_{k, L}$ of holomorphic modular forms. On the space $V_{k,L} \times V_{2-k,L^-}$ we define the bilinear pairing 
\begin{equation}
  \label{eq:pairing}
    \{g, f \} :=  \CT_{q = 0} \sum_{h \in A_L} g_h f_h. 
\end{equation}
When restricted to $M^!_{k, L}  \times    M^!_{2-k, L^-}$, this pairing vanishes identically. 
Furthermore, if $g \in M^!_{k,L}$, and $f = \widetilde{G}^+$ is the holomorphic part of a harmonic Maass form $\widetilde{G}$ of weight $2-k$ for $\rho_{L^-}$ such that $\xi_k \widetilde{G} = G \in M_{k, L}$, an application of Stokes' theorem gives us (see e.g.\ Proposition~3.5 in \cite{BF04})
\begin{equation}
  \label{eq:Stokes}
  \{g , \widetilde{G}^+ \}=   (g, G)^\mathrm{reg}
  := \lim_{T \to \infty} \int_{\mathcal{F}_T} \langle g(\tau), G(\tau) \rangle  v^k \frac{dudv}{v^2},
\end{equation}
where $\mathcal{F}_T$ is the usual fundamental domain for $\SL_2(\Zb)$ truncated at the height $T$. It is clear that the pairing vanishes on $V^+_{k,L} \times S_{2-k, L^-}$.
Denote $W_{k, L} \subset V_{k, L}$ the subspace spanned by $V^+_{k, L}$ and $M^!_{k, L}$.
By Serre duality (compare Theorem 3.1 in \cite{Borcherds99}), we know that 
\begin{equation}
  \label{eq:Serre}
  S_{2-k, L^-}  = W_{k, L}^\perp \subset V_{2-k, L^-},\qquad  
  W_{k, L} = S_{2-k, L^-}^\perp \subset V_{k,L},
\end{equation}
where the orthogonal complement is taken with respect to the pairing $\{\cdot,\cdot\}$. From this, we can deduce the following result.

\begin{lemma}
  \label{lemma:dual}
  The pairing $\{\cdot, \cdot\}$ induces a perfect pairing on
$    V_{k, L} / M^!_{k, L}  \times    M^!_{2-k, L^-}$.
\end{lemma}

\begin{proof}
Suppose $f \in (M^!_{2-k, L^-})^\perp \subset V_{k, L}$. Then $f \in S_{2-k, L^-}^\perp = W_{k, L}$ by \eqref{eq:Serre}, and we can write $ f= f_1 + f_2$ with $f_1\in V^+_{k, L}$ and $f_2 \in M^!_{k, L} \subset (M^!_{2-k, L^-})^\perp$. 
That means $f_1$ is in $(M^!_{2-k, L^-})^\perp \cap (V_{2-k, L^-}^+)^\perp$, and hence in $W_{2-k, L^-}^\perp$. 
Again by \eqref{eq:Serre}, we have $f_1 \in S_{k, L}$, hence $f = f_1 + f_2 \in M^!_{k, L}$. 
\end{proof}

\begin{prop}
  \label{prop:xipreimage}
  Let $G \in M_{k, L}$, and suppose that $\widetilde{G}^+ \in V_{2-k, L^-}$ satisfies
 \begin{equation}
    \label{eq:pairing2}
    \{ g, { \widetilde{G}^+} \}
    = ( g, G )^{\reg}    
  \end{equation}
  for all $g  \in M^!_{k, L}$. Then $\widetilde{G}^+$ is a mock modular form with shadow $G$. 
\end{prop}

\begin{proof}
  The case $G = 0$ follows directly from Lemma \ref{lemma:dual}.
  More generally, we can subtract from $\widetilde{G}^+$ a known mock modular form with shadow $G$ to reduce to the case $G = 0$. 
\end{proof}

\subsection{Hecke-type Operators}
\label{subsec:Hecke-op}
We quickly recall certain Hecke-type operators on vector-valued modular forms (see e.g.\ \cite[section 2]{BS20}).
Given $k \in \Zb + \frac{1}{2}$, $N, d \in \Nb$ and $\e = \pm 1$, there are operators $U_d : H_{k, \Zb[N]^\e} \to H_{k, \Zb[Nd^2]^\e}$ and $V_d : H_{k, \Zb[N]^\e} \to H_{k, \Zb[Nd]^\e}$ such that 
for $f(\tau) = f^+(\tau) + f^-(\tau) \in H_{k, \Zb[N]^\e}$ with $f^+(\tau) = \sum_{m, h} c(h, m) q^m \ef_h$, we have
\begin{equation}
  \label{eq:Hecke-op-U}
  \begin{split}
    (    f \mid U_d)^+ &= \sum_{h \in A_{Nd^2}} \ef_h \sum_{m \in \Qb} c(m/d^2, h/d) q^{m/d^2},\\
    \xi_k(f \mid U_d) &= (\xi_k f) \mid U_d,
  \end{split}
\end{equation}
\begin{equation}
  \label{eq:Hecke-op-V}
  \begin{split}
    (    f \mid V_d)^+ &= \sum_{h \in A_{Nd}} \ef_h \sum_{m \in \Qb} \sum_{a \mid \lp (m - \e\tfrac{h^2}{4N})/p, h, d\rp} a^{k-1/2}  c(m/a^2, h/a) q^{m/d},\\
    \xi_k(f \mid V_d) &= d^{k-1} (\xi_k f) \mid V_d.
  \end{split}
\end{equation}
These two operators commute with each other for all $d$.
When we identify holomorphic scalar-valued modular forms for $\rho_N$ with holomorphic Jacobi forms of index $N$, the operators $U_d, V_d$ are the usual index changing operators (see \cite[section I.4]{EZ85}).
For any $c \parallel N$, the automorphism $W_c \in \mathrm{Aut}(A_N)$ defined by
\begin{equation}
  \label{eq:sigmac}
  W_c(h) \equiv - h \bmod 2c,~
  W_c(h) \equiv  h \bmod (N/2c),
\end{equation}
acts on $H_{k, \Zb[N]^\e}$ through its components, which commutes with $\xi_k$. Clearly, the three operators $U_d, V_d, W_c$ all preserve integrality of the holomorphic part Fourier coefficients.
As in \cite[section 4.4]{DMZ12}, we also define $  \Vc_{d}^{(N)}: H_{k, \Zb[N]^\e} \to H_{k, \Zb[Nd]^\e}$ by
\begin{equation}
  \label{eq:Vc}
\Vc_{d}^{(m)} := \sum_{r^2 \mid d,~ (r, m) = 1} \mu(r) r^{k-1/2} V_{d/r^2} U_r,
\end{equation}
where $\mu(\cdot)$  is the M\"obius function.
We then have the following result.
\begin{lemma}
  \label{lemma:Hecke-op}
  For each $N \in \Nb$ and $k \in \{\tfrac{1}{2}, \tfrac{3}{2}\}$, let $\phi_{k, N}(\tau) \in H_{k, \Zb[N]^-}$  be harmonic Maass forms such that
  $\xi_{k} \phi_{k, N}(\tau)  = N^{k-1}   \theta_N(\tau; 3/2 - k)$.
For $N_0 \in \Nb$ and prime $p \nmid N_0$, the modular form
  $$
  \phi_{k, N_0 p^e } - \lp \phi_{k, N_0} \mid \Vc_{p^e}^{(1)} + \phi_{N_0 p} \mid (1 - W_p) \Vc_{p^{e-1}}^{(p)}
  \rp/2
  $$
is weakly holomorphic for any $e \ge 1$ and $k \in \{\tfrac{1}{2}, \tfrac{3}{2}\}$.
\end{lemma}

\begin{proof}
  This is simply equation (10.48) in \cite{DMZ12} phrased in terms of vector-valued modular forms, and follows from the same calculations loc.\ cit.
\end{proof}

\section{Theta lifts}
\label{sec:theta lifts}

In this section, we compute some regularized theta lifts for lattices of signature $(1, n)$ for $n \geq 1$. The formulas in this section are special cases of the general results of Borcherds \cite{Borcherds98} and Bruinier \cite{bruinierhabil}, but we write down the simplifications for the convenience of the reader.

Throughout this section, we let $L$ be an even lattice of signature $(1,n)$ with $n \geq 1$, and we fix
an isotropic vector $\ell \in V_\Rb:= L \otimes \Rb$. If $L$ is isotropic, we will choose $\ell$ to be a primitive isotropic vector in $L$. We let $\Gr(L)$ denote the Grassmannian of positive lines in $V_\Rb$. For $z \in \Gr(L)$ we consider the polynomials
 \[
 p_z(\lambda) := \left(\lambda,\frac{\ell_z}{|\ell_z|}\right), \qquad p_{z^\perp}(\lambda) := \left(\lambda,\frac{\ell_{z^\perp}}{|\ell_{z^\perp}|}\right),
 \] 
 on $V_\Rb$, where we write $\lambda_z$ for the projection of $\lambda \in V_\Rb$ to the subspace $z$, and $|\lambda| := \sqrt{|(\lambda,\lambda)|}$. For $m^+, m^- \in \{0,1\}$ we define the Siegel theta function
\begin{align}\label{eq:theta function}
\Theta_{L, \ell}^{(m^+,m^-)}(\tau,z) := v^{\frac{n}{2}+m^-}\sum_{\lambda \in L'}
p_z^{m^+}(\lambda)p_{z^\perp}^{m^-}(\lambda)
\ebf\left(Q(\lambda_z)\tau + Q(\lambda_{z^\perp})\overline{\tau}\right)\ef_{\lambda+L}
\end{align}
on $\Hb \times \Gr(L)$. Note that for $m^+,m^- \in \{0,1\}$ and every fixed $z \in \Gr(L)$ the polynomial $p_z^{m^+}(\lambda)p_{z^\perp}^{m^-}(\lambda)$ is harmonic and homogeneous of degree $(m^+,m^-)$. Hence by Theorem~4.1 in \cite{Borcherds98} the theta function transforms like a modular form of weight $\frac{1-n}{2}+m^+-m^-$ for $\rho_L$ in $\tau$. Moreover, it is $\Gamma_L$-invariant in $z$. The corresponding regularized theta lift of a weakly holomorphic modular form $f \in M_{\frac{1-n}{2}+m^+-m^-,L}^!$ is defined by
 \begin{equation}
  \label{eq:theta lift}
\Phi_{L, \ell}^{(m^+,m^-)}(f,z) := \lim_{T \to \infty}\int_{\mathcal{F}_T}\left\langle f(\tau),\Theta_{L, \ell}^{(m^+,m^-)}(\tau,z) \right\rangle 
v^{\frac{1-n}{2}+m^+-m^-}\frac{du dv}{v^2}.
\end{equation}
Note that when $m^- = 0$ (resp.\ $m^- = 1 = n$), we will fix generators of $z$ (resp.\ $z$ and $z^\perp$) to remove the dependence on $\ell$, and omit it from the subscripts.

By the general theory developed in \cite{Borcherds98,bruinierhabil}, the theta lift converges for every $z \in \Gr(L)$ and is real analytic up to singularities along the Heegner divisor
\[
Z_f := \sum_{\substack{\lambda \in L' \\ Q(\lambda) < 0}}a_f(\lambda,Q(\lambda))\lambda^\perp,
\] 
where $\lambda^\perp$ denotes the hypersurface consisting of all $z \in \Gr(L)$ perpendicular to $\lambda$. 
These hypersurfaces partition $\Gr(L)$ into infinitely many connected components, the so-called Weyl chambers corresponding to $f$.

\subsection{Theta lifts on isotropic lattices of signature $(1,n)$}\label{theta lifts}

Let $L$ be an isotropic even lattice of signature $(1,n)$ with $n \geq 1$. Let $\ell \in L$ be primitive isotropic and let $M \in \Nb$ with $(\ell,L) = M\Zb$. Choose some $\ell' \in L'$ with $(\ell,\ell') = 1$, some $\zeta \in L$ with $(\ell,\zeta) = M$, and set
\[
K := L \cap \ell^\perp \cap \ell'^\perp.
\]
Then $K$ has signature $(0,n-1)$ and $L = K \oplus \Zb \zeta \oplus \Zb \ell$ (see Proposition~2.2 in \cite{bruinierhabil}). Let 
\[
L_0' := \{\lambda \in L': (\lambda,\ell) \equiv 0 \pmod{M}\},
\]
and let $p: L_0' \to K'$ be defined by
\begin{equation}
  \label{eq:p}
  p(\lambda) := \lambda_K-\frac{(\lambda,\ell)}{M}\zeta_K, \qquad
  \lambda_K := \lambda + ((\lambda, \ell)(\ell', \ell') - (\lambda , \ell')) \ell - (\lambda, \ell) \ell' \in K \otimes \Qb,
\end{equation}
which induces a surjection $p: L_0'/L \to K'/K$. 
For $h \in L_0'/L$, routine calculations give us
\begin{equation}
  \label{eq:ph}
p(h)     = h - (h , \ell') \ell \in V/L.  
\end{equation}

We let $\mathbb{B}_k(x)$ denote the one-periodic function that agrees with the Bernoulli polynomial $B_k(x)$ for $0 \leq x < 1$. Recall that the first few Bernoulli polynomials are given by
\[
B_0(x) = 1, \qquad B_1(x) = x-\frac{1}{2}, \qquad B_2(x) = x^2-x+\frac{1}{6}, \qquad B_3(x) = x^3 - \frac{3}{2}x^2 + \frac{1}{2}x.
\]
The theta lift has the following Fourier expansion.

\begin{prop}\label{prop:fourier expansion isotropic}
	Let $m^+,m^- \in \{0,1\}$ and let $L$ be an even lattice of signature $(1,n)$ with a primitive isotropic vector $\ell \in L$. For $f \in M_{\frac{1-n}{2}+m^+-m^-,L}^!$ and every $z \in \Gr(L)$ we have
	\begin{align*}
	\Phi_{L, \ell}^{(m^+,m^-)}(f,z) 
          &= \frac{|\ell_z|^{m^+-1}(-|\ell_{z^\perp}|)^{-m^-}}{\sqrt{2}(4\pi)^{m^+}}\Phi_K^{(m^+,m^-)}(f)
  -\frac{\sqrt{2}(-4\pi)^{1-m^+}}{2-m^++m^-}|\ell_z|^{1-m^++2m^-}(-|\ell_{z^\perp}|)^{-m^-} \\
	&\qquad \times \sum_{\lambda \in K'}\sum_{\substack{\delta \in L_0'/L \\ p(\delta) = \lambda + K}}a_f(\delta,Q(\lambda))\mathbb{B}_{2-m^++m^-}\left(\frac{p_z(\lambda)}{|\ell_z|}+(\delta,\ell')\right) \\
	&\quad + \frac{(-4\pi)^{1-m^+}}{\sqrt{2}}\sum_{\substack{\lambda \in L' \\ Q(\lambda) < 0 \\ (\lambda,\ell) \neq 0}}
a_f(\lambda,Q(\lambda))p_z^{1-m^+}(\lambda)p_{z^\perp}^{m^-}(\lambda)\big(\sgn(p_z(\lambda)) - \sgn((\lambda, \ell)) \big),
	\end{align*}
	where the constant $\Phi_K^{(m^+,m^-)}(f)$ is given by
	\begin{align*}
	\Phi_K^{(m^+,m^-)}(f) = \begin{dcases}
	-8\pi\sum_{m\geq 0}\sum_{\lambda \in K'}\sum_{\substack{\delta \in L_0'/L \\ p(\delta)=\lambda+K}}a_f(\delta,-m+Q(\lambda))\sigma_1(m) & \text{if } m^+ = m^-, \\
	0 & \text{otherwise},
	\end{dcases}
	\end{align*}
	with $\sigma_1(m) :=\sum_{d \mid m}d$ for $m \in \Nb$ and $\sigma_1(0) := -\frac{1}{24}$.
\end{prop}

\begin{rmk}\label{remark finite singular part}
	By the same arguments as in the proof of Theorem 3.3 in \cite{BS19} one can show that the sum in the third line in Proposition~\ref{prop:fourier expansion isotropic} is finite for every fixed $z \in \Gr(L)$, and vanishes for $|\ell_z|$ small enough. Moreover, the second line encompasses the singularities of the theta lift along those $\lambda^\perp$ with $(\lambda,\ell) \neq 0$, whereas the first line gives those with $(\lambda,\ell) = 0$.
\end{rmk}

\begin{rmk}
	For $m^+ = m^-$ and $n \geq 2$, the constant $\Phi_K^{(m^+,m^-)}(f)$ can also be computed using Theorem~2.14 in \cite{bruinierhabil} by writing $f$ as a linear combination of the Maass Poincar\'e series defined in \eqref{eq:maass poincare series}. This yields the alternative representation 
	\begin{align}\label{eq alternative first line}
	\Phi_K^{(m^+,m^-)}(f) = \frac{8\pi}{(n-1)}\sum_{\lambda \in K'}\sum_{\substack{\delta \in L_0'/L \\ p(\delta)=\lambda+K}}a_f(\delta,Q(\lambda))|Q(\lambda)|.
	\end{align}
\end{rmk}

\begin{rmk}
  \label{rmk:simplify} 
  For any generator $w \in z$ and $\tilde\lambda \in L + \delta$ such that $p(\tilde\lambda) = \lambda$, we have
  \begin{equation}
    \label{eq:Bsim}
    \frac{p_z(\lambda)}{|\ell_z|}+(\delta,\ell') =
  \frac{( \tilde\lambda, w)}{(\ell, w)}
    - \frac{(\tilde\lambda,\ell)}{M} \cdot \frac{(\zeta, w)}{(\ell, w)} \in \Rb/\Zb.  
  \end{equation}
\end{rmk}
\begin{rmk}
  \label{rmk:Flz}
  There is $F_{\ell, z} \in V_{\frac{3+n}{2} - m^+ + m^-, L^-}$ that allows us to rewrite Proposition~\ref{prop:fourier expansion isotropic} as
  $$
  \Phi^{(m^+, m^-)}_{L, \ell} = \left\{
f, F_{\ell, z}
    \right\}
    $$
    for all $f \in M^!_{\frac{1-n}{2} + m^+ - m^-,L}$ (see e.g.\ \eqref{eq:Flw} for $n = 4, m^+ = 1, m^- = 0$).
\end{rmk}
\begin{proof}[Proof of Proposition~\ref{prop:fourier expansion isotropic}]
	The formula follows from Theorem~10.2\footnote{Note that there is a sign $(-1)^{h^-}$ missing in the cited formula.} in \cite{Borcherds98} (or a generalization of Proposition~3.1 in \cite{bruinierhabil}), where the Fourier expansion of $\Phi_{L, \ell}^{(m^+,m^-)}(f,z)$ was computed in a fixed Weyl chamber of $\Gr(L)$. In order to extend the expansion to all of $\Gr(L)$, one can use the shape of the singularities of the theta lift to determine its ``wall crossing'' behavior as a point $z \in \Gr(L)$ moves across a hypersurface in the Heegner divisor $Z_f$ from one Weyl chamber to another (compare also Corollary~6.3 in \cite{Borcherds98}). 
	
	First, by Theorem~10.2 in \cite{Borcherds98}, for $|\ell_z|$ small enough the Fourier expansion of the theta lift is given by the first three lines of the expression on the right-hand side in the proposition. The constant $\Phi_K^{(m^+,m^-)}(f)$ appearing in Theorem~10.2 in \cite{Borcherds98} vanishes if $m^+ \neq m^-$, and for $m^+ = m^-$ it is equal to the regularized integral
	\begin{align}\label{eq average value}
	\Phi_K^{(m^+,m^-)}(f) = \lim_{T \to \infty}\int_{\mathcal{F}_T}\left\langle f_K(\tau), \overline{\Theta_{K^-}(\tau)}\right\rangle \frac{du dv}{v^2}, \qquad f_{K}(\tau) := \sum_{\gamma \in K'/K}\sum_{\substack{\delta \in L_0'/L \\ p(\delta) = \gamma}}f_{\delta}(\tau)\ef_{\gamma},
	\end{align}
	where $\Theta_{K^-}(\tau) = \sum_{\lambda \in K'}\ebf(-Q(\lambda)\tau)\ef_{\lambda+K}$ is the holomorphic theta function associated to the positive definite lattice $K^-$. In particular, the integral in \eqref{eq average value} can be viewed as the regularized average value of a weakly holomorphic modular form of weight $0$ for $\SL_2(\Zb)$. Hence \eqref{eq average value} can be evaluated as explained in Remark~4.9 in \cite{BF06}.
	
	Moreover, by Theorem 6.2 in \cite{Borcherds98} (or a generalization of Theorem~2.12 in \cite{bruinierhabil}), the theta lift has a singularity of type
	\[
	\frac{(-4\pi)^{1-m^+}}{\sqrt{2}}\sum_{\substack{\lambda \in L' \\ Q(\lambda) < 0 \\ \lambda \perp z_0}}
a_f(\lambda,Q(\lambda))p_z^{1-m^+}(\lambda)p_{z^\perp}^{m^-}(\lambda)\sgn(p_z(\lambda))
	\]
	at a point $z_0 \in \Gr(L)$. Here we say that a function $f$ has a singularity of type $g$ at a point $z_0$ if there exists a neighborhood $U$ of $z_0$ such that $f$ and $g$ are defined on $U$ and $f-g$ is real analytic on $U$. This definition slightly differs from the one used in \cite{Borcherds98} and \cite{bruinierhabil}, but it allows us to extend the Fourier expansion to points $z \in \Gr(L)$ where the theta lift has singularities.
	
	It is now easy to check that the expression on the right-hand side of the proposition has the same singularities as $\Phi_{L, \ell}^{(m^+,m^-)}(f,z)$. In particular, the difference of $\Phi_{L, \ell}^{(m^+,m^-)}(f,z)$ and the expression in the proposition defines a real-analytic function on all of $\Gr(L)$ and vanishes for $|\ell_z|$ small enough (see also Remark~\ref{remark finite singular part}), and hence vanishes everywhere on $\Gr(L)$. This finishes the proof.
\end{proof}

\subsection{Theta lifts on anisotropic lattices of signature $(1,1)$}

We now compute the theta lift $\Phi_{L, \ell}^{(m^+,m^-)}(f,z)$ defined in~\eqref{eq:theta lift} for anisotropic lattices $L$ of signature $(1,1)$. One can compute the theta lift for anisotropic lattices of signature $(1,n)$ for any $n \geq 1$ in a similar way, but the resulting formulas do not look as pleasing. Hence we confine ourselves with signature $(1,1)$, which suffices for our applications.
Another advantage is that we can fix generators of $z, z^\perp$ and remove the isotropic vector $\ell$ from the notations.

\begin{prop}\label{prop:unfolding anisotropic}
	Let $m^+,m^- \in \{0,1\}$ but $(m^+,m^-) \neq (1,0)$ and let $L$ be an anisotropic lattice of signature $(1,1)$. 
For any $f \in M_{m^+-m^-,L}^!$ and $z \in \Gr(L)$ we have
	\begin{align*}
	&\Phi_{L}^{(m^+,m^-)}(f,z) \\
	&\quad = \frac{2^{\frac{3}{2}-m^+-m^-}}{\pi^{m^+-1}}\sum_{\substack{\lambda \in L' \\ Q(\lambda) < 0}}a_f(\lambda,Q(\lambda))\sgn(p_z(\lambda))^{m^+}\sgn(p_{z^{\perp}}(\lambda))^{m^-}(|p_{z^\perp}(\lambda)|-|p_z(\lambda)|)^{1-m^++m^-}.
	\end{align*}
\end{prop}

\begin{rmk}
	For $(m^+,m^-) = (1,0)$ the proof below does not work since there might be non-trivial holomorphic modular forms of weight $1$ for $\rho_L$, so we cannot write $f$ as a linear combination of Maass Poincar\'e series (and, possibly, invariant vectors). Hence we exclude this case from the above proposition.
\end{rmk}

\begin{proof}[Proof of Proposition~\ref{prop:unfolding anisotropic}]	
	First note that we cannot use Theorem~10.2 in \cite{Borcherds98} since it requires the existence of an isotropic vector in $L$. Instead, we will write
	\[
	f(\tau) = \frac{1}{2}\sum_{h \in A_L}\sum_{m < 0}a_f(h,m)F_{h,m}(\tau) + \mathfrak{c}
	\]
	as a linear combination of the Maass Poincar\'e series $F_{h,m}$ defined in \eqref{eq:maass poincare series} and (if $m^+ = m^-$) a $\rho_L$-invariant vector $\mathfrak{c} \in \Cb[A_L]$. Then we compute the lift of $F_{h,m}$ and $\mathfrak{c}$ using the unfolding argument as in the proof of Theorem~2.14 in \cite{bruinierhabil}.
	
	To simplify the notation, we only treat the case $m^+ = m^- = 1$ here. The other cases are analogous. We first show that the theta lift $\Phi_{L}^{(1,1)}(\mathfrak{c},z)$ of an invariant vector $\mathfrak{c} \in \Cb[A_L]$ vanishes identically. To this end, we use the simple fact that every invariant vector can be written as a linear combination of residues at $s = 1$ of Eisenstein series
	\[
	E_{h}(\tau,s) := \sum_{(A,\phi) \in \widetilde{\Gamma}_\infty \backslash \Mp_2(\Zb)}v^s \ef_h |_{0,L}(A,\phi)
	\]
	corresponding to isotropic elements $h \in A_L$. By the usual unfolding argument one can show that $\Phi_{L}^{(1,1)}(E_h(\cdot,s),z)$ is a multiple of
	\[
	\sum_{\substack{\lambda \in (L+h) \setminus \{0\} \\ Q(\lambda) = 0}}\frac{p_z(\lambda)p_{z^\perp}(\lambda)}{(4\pi|Q(\lambda_{z^\perp})|)^{s+\frac{1}{2}}}
	\]
	for $\Re(s)$ large enough. Since $L$ is anisotropic, the sum over $\lambda \in (L+h) \setminus \{0\}$ with $Q(\lambda) = 0$ is empty, so the theta lift of $E_h(\tau,s)$ vanishes identically for $\Re(s)$ big enough. In particular, its residue at $s = 1$ vanishes, as well. This shows that the theta lift of an invariant vector vanishes. 
	
	Hence it suffices to compute the lifts of Poincar\'e series $F_{h,m}$. Using the unfolding argument again, we find for $\Re(s)$ big enough
	\begin{align*}
	&\Phi_{L}^{(1,1)}(F_{h,m}(\cdot,s),z) \\
	&= \frac{2}{\Gamma(2s)}\sum_{\substack{\lambda \in L+h \\ Q(\lambda) = m}}p_z(\lambda)p_{z^{\perp}}(\lambda)\int_{0}^{\infty}v^{-\frac{1}{2}}M_{0,s-\frac{1}{2}}(4\pi |m| v)\exp\lp-2\pi v(Q(\lambda_z)-Q(\lambda_{z^\perp}))\rp dv.
	\end{align*}
	The integral is an inverse Laplace transform (see equation (11) on p. 215 of \cite{tables}) given by
	\begin{align*}
	\frac{(4\pi |m|)^s \Gamma\lp s+\frac{1}{2}\rp}{(4\pi |Q(\lambda_{z^\perp})|)^{s+\frac{1}{2}}}\ _2 F_{1}\left(s+\frac{1}{2},s,2s;\frac{|m|}{|Q(\lambda_{z^\perp})|} \right).
	\end{align*}
	Plugging in $s = 1$ and using that $_2 F_{1}\left(\frac{3}{2},1,2;x\right) = 2\frac{1-\sqrt{1-x}}{x\sqrt{1-x}}$ we find
	\[
	\Phi_{L}^{(1,1)}(F_{h,m},z) = \sum_{\substack{\lambda \in L+h \\ Q(\lambda) = m}}p_z(\lambda)p_{z^\perp}(\lambda)\frac{\sqrt{|Q(\lambda_{z^\perp})|}-\sqrt{Q(\lambda_{z})}}{\sqrt{|Q(\lambda_{z^\perp})|}\sqrt{Q(\lambda_{z})}}.
	\]
	Using $Q(\lambda_z) = \frac{1}{2}p_z(\lambda)^2$ and $Q(\lambda_{z^\perp}) = -\frac{1}{2}p_{z^\perp}(\lambda)^2$ we obtain the stated formula.
	\end{proof}
	
	Note that, in contrast to the isotropic case, the sum on the right-hand side of Proposition~\ref{prop:unfolding anisotropic} is not finite since the discriminant kernel $\Gamma_L$ is infinite (as it corresponds to a non-trivial subgroup of the units in a real quadratic field). However, one can obtain a finite evaluation of the theta lift by splitting the sum over $\lambda \in L'$ modulo $\Gamma_L$. To this end, it is convenient to view anisotropic lattices of signature $(1,1)$ as lattices in real quadratic fields, as we now explain.
	
	Let $D > 1$ be a non-square discriminant (not necessarily fundamental), let $F = \Qb(\sqrt{D})$ be the corresponding real quadratic field, and let $\mathcal{O}_F$ be its ring of integers. We consider the subring of $\mathcal{O}_F$ given by 
\[
\mathcal{O}_D := \Zb + \frac{D+\sqrt{D}}{2}\Zb.
\]
Notice that $\mathcal{O}_D = \mathcal{O}_F$ if $D$ is a fundamental discriminant. More generally, if $D = f^2 D_F$ with a fundamental discriminant $D_F$ and $f \in \Nb$, then $\mathcal{O}_D$ is the order of discriminant $D$ and conductor $f$ in $F$. For an integral ideal $\mathfrak{a} \subseteq \mathcal{O}_D$ with $A := [\mathcal{O}_D:\mathfrak{a}]$ and a positive integer $M \in \Nb$ we consider the lattice
\begin{equation}
  \label{eq:LaM}
(L_{\mathfrak{a},M},Q_{\mathfrak{a},M}) := \left(M\mathfrak{a},-\frac{\Nm_{F/\Qb}}{AM} \right),
\end{equation}
where $\Nm$ denotes the norm. The associated bilinear form is $B_{\mathfrak{a},M}(\lambda,\mu) = -\frac{\Tr_{F/\Qb}(\lambda \cdot \mu')}{AM}$, where $\mu'$ denotes the conjugate of $\mu$. The lattice $L_{\mathfrak{a},M}$ is anisotropic of signature $(1,1)$, and contains the sublattice
\begin{equation}
  \label{eq:iota}
  \begin{split}
\iota:    \Zb[AMD] \oplus \Zb[AM]^- &\hookrightarrow L_{\af, M}\\
    [b, a] &\mapsto AM(a + b\sqrt{D}).
  \end{split}
\end{equation}
The dual lattice is given by $\mathfrak{a}\mathfrak{d}_D^{-1}$, where $\mathfrak{d}_D := \sqrt{D}\mathcal{O}_D$, so the discriminant group of $L_{\mathfrak{a},M}$ is isomorphic to $\mathcal{O}_D/M\mathfrak{d}_D$. The discriminant kernel $\Gamma_{L_{\mathfrak{a},M}}$ is generated by $-1$ and a totally positive unit $\varepsilon_{L_{\mathfrak{a},M}} > 1$ in $\mathcal{O}_D$, which is an integral power of the smallest totally positive unit $> 1$ in $\mathcal{O}_D$.

	For simplicity, in the following corollary we already evaluate the theta lift at a certain special point which is adjusted to our later applications. We view $(L_{\mathfrak{a},M},Q_{\mathfrak{a},M})$ as a sublattice of $(\Rb^2,xy)$ by sending $\lambda \in L_{\mathfrak{a},M}$ to $\frac{1}{\sqrt{AM}}[\lambda,-\lambda']$. We will choose the generators $(z_1,z_2),(z_1^\perp,z_2^\perp) \in \Rb^2$ of $z$ and $z^\perp$ with $z_1,z_2 > 0$ and $z_1^\perp > 0, z_2^\perp < 0$.
	
	\begin{cor}\label{cor:unfolding anisotropic}
	Suppose that $(L,Q) = (L_{\a,M},Q_{\a,M})$ is an anisotropic lattice of the form \eqref{eq:LaM}. Let $f \in M_{m^+-m^-,L}^!$ and let $w \in \Gr(L)$ be the positive line generated by $\frac{1}{\sqrt{2}}[1,1]$. Then we have
	\begin{align*}
	\Phi_{L}^{(1,1)}(f,w)  &= -\frac{1}{\sqrt{AM}}\sum_{\substack{\lambda \in R_L}}a_f(\lambda,Q(\lambda))\sgn(\lambda)\Tr_{F/\Qb}\left(\frac{\lambda}{1-\varepsilon_L^{-1}} \right), \\
	\Phi_{L}^{(0,0)}(f,w)  &= \frac{4 \pi}{\sqrt{DAM}}\sum_{\substack{\lambda \in R_L^*}}a_f(\lambda,Q(\lambda))\sgn(\lambda)\Tr_{F/\Qb}\left( \frac{\sqrt{D}\lambda}{1-\varepsilon_L^{-1}}\right), \\
	\Phi_{L}^{(0,1)}(f,w)  &= \frac{2\sqrt{2} \pi}{\sqrt{D}AM}\sum_{\substack{\lambda \in R_L^*}}a_f(\lambda,Q(\lambda))\sgn(\lambda)\Tr_{F/\Qb}\left(\frac{\sqrt{D} \lambda^2}{1-\varepsilon_L^{-2}}\right),
	\end{align*}
	where $R_L := \{\lambda \in L': \varepsilon_L^{-2} < \lambda/\lambda' < 1\}$ and $R_L^* := \{\lambda \in L': \varepsilon_L^{-2} < \lambda/\lambda' \leq 1\}$. The sums on the right-hand sides are finite.
	\end{cor}
	
	\begin{proof}
	Again, we only treat the case $(m^+,m^-) = (1,1)$ for simplicity. First note that we have
	\[
	p_w(\lambda) = \frac{1}{\sqrt{2AM}}(\lambda-\lambda'),\qquad p_{w^\perp}(\lambda) = \frac{1}{\sqrt{2AM}}(\lambda+\lambda').
	\]
	Using Proposition~\ref{prop:unfolding anisotropic} we find
	\begin{align*}
	\Phi_{L}^{(1,1)}(f,w) &= \frac{1}{2\sqrt{AM}}\sum_{\substack{\lambda \in L' \\ Q(\lambda) < 0}}a_f(\lambda,Q(\lambda))\sgn\lp |\lambda'| - |\lambda|\rp \big(|\lambda+\lambda'|-|\lambda-\lambda'| \big) \\
	&= -\frac{1}{\sqrt{AM}}\sum_{\substack{\lambda \in L' \\ Q(\lambda) < 0}}a_f(\lambda,Q(\lambda))\sgn\lp |\lambda'| - |\lambda|\rp \min(|\lambda|,|\lambda'|),
	\end{align*}
	where we used that $|x+y|-|x-y| = -2\min(|x|,|y|)$ for $x,y \in \Rb$ with $xy < 0$, and
	$\sgn(\lambda) = \sgn(\lambda')$ since $Q(\lambda) = -\frac{\lambda\lambda'}{AM} < 0$. Let $\Gamma_L'$ be the subgroup of $\Gamma_L$ consisting of totally positive units.
Note that the terms with $\lambda = \lambda'$ contribute nothing.
 As a system of representatives for $\Gamma_L' \backslash\{\lambda \in L': Q(\lambda) < 0, \lambda \neq \lambda'\}$ we choose the set $R_L$. Then the set of all $\lambda \in L'$ with $Q(\lambda) < 0$ and $\lambda \neq \lambda'$ is given by $\{\lambda \varepsilon_L^n: \lambda \in R_L\}$. For $\lambda \in R_L$ and $n \in \Zb$ we have
	\[
	\min(|\lambda\varepsilon_L^n|, |\lambda' \varepsilon_L^{-n}|) = \begin{cases}
	|\lambda| \varepsilon_L^n & \text{if } n\leq 0, \\
	|\lambda'|\varepsilon_L^{-n}& \text{if } n > 0,
	\end{cases}
	\qquad \sgn(|\lambda' \varepsilon_L^{-n}|-|\lambda \varepsilon_L^{n}|)= \begin{cases}
	1 & \text{if } n\leq 0, \\
	-1& \text{if } n \geq 1.
	\end{cases}
	\]
	We obtain
	\begin{align*}
	\sum_{\substack{\lambda \in L' \\ Q(\lambda) < 0}}a_f(\lambda,Q(\lambda))\sgn\lp |\lambda'| - |\lambda|\rp \min(|\lambda|,|\lambda'|) &= \sum_{\lambda \in R_L}a_f(\lambda,Q(\lambda))\left(|\lambda| \sum_{n\leq 0}\varepsilon_L^n - |\lambda'|\sum_{n \geq 1}\varepsilon_L^{-n}
 \right) \\
	&= \sum_{\lambda \in R_L}a_f(\lambda,Q(\lambda))\left(\frac{|\lambda|}{1-\varepsilon_L'}+\frac{|\lambda'|}{1-\varepsilon_L}
\right).
	\end{align*}
	This yields the stated formula. The sums on the right-hand side of the corollary are finite since $f$ has finite principal part and the intersection of $R_L$ (resp. $R_L^*$) with the set of vectors of a fixed norm is finite.
\end{proof}

\section{Proof of Theorem~\ref{thm:denom}}
\label{sec:proof theorem denom}

We will give the full proof for case of weight $k = 1/2$ harmonic Maass form.
The weight $3/2$ case follows from the same argument.

For $N \in \Nb$, define 
  \begin{equation}
    \label{eq:db}
    \begin{split}
    S(N)&:=\{    \phi \in H_{1/2, \Zb[N]^-}: \xi_{1/2}(\phi) =\tfrac{1}{\sqrt{N}}\theta_N(\tau; 1)\},\\
    \db(N) &:= \min\{m \in \Nb: m \phi^+ \text{ has integral Fourier coefficients for some } \phi \in S(N)\}.
    \end{split}
  \end{equation}
  Remark \ref{rmk:denom} implies that we can replace min in the definition of $\db(N)$ with gcd.
From the construction in equation (9.2) of \cite{DMZ12}, we know that 
\begin{equation}
  \label{eq:DMZbound}
  \db(N) \mid 2 \cdot 12^{N-1}.
\end{equation}
for all $N \in \Nb$.

  Consider the even lattice $L = \Zb[N] \oplus (\Zb[3]^-)^4$ of signature $(1,4)$ and level $12N$, and let $N_d$ be the same as in Theorem \ref{thm:denom}.
By the four-square theorem, we can find a primitive isotropic vector $\ell = [3/N_3,*] \in L$. Then we have $(\ell, L) = M\Zb$ for some $M \in \Nb$ satisfying $M \mid 6N$.
  Let $\zeta \in L$ and $\ell' = [\ell'_1, *] \in L'$  such that $(\zeta, \ell) = M$ and $(\ell, \ell') = 1$.
  Denote $K = L \cap \ell^\perp \cap \ell'^\perp$.

  Now we come to the construction of a harmonic Maass form $\phi \in S(N)$.
  To this end, we will first compute the regularized Petersson inner product $(g(\tau),\theta_N(\tau;1))^\reg$ for every $g \in M_{\frac{3}{2},\Zb[N]}^!$, using the theta lift studied in Section~\ref{sec:theta lifts}. Then we can express this inner product in terms of the pairing $\{g,\widetilde{\theta}_{N}^+(\tau;1)\}$ for some explicit Laurent series $\widetilde{\theta}_{N}^+(\tau;1) \in V_{\frac{1}{2},\Zb[N]^-}$, and finally apply Proposition~\ref{prop:xipreimage} to obtain the desired mock modularity.
   
   It is easy to check that the theta function $ \Theta_{L, \ell}^{(1,0)}(\tau,z)$ defined in \eqref{eq:theta function} splits at the line generated by 
   \[
   w = \frac{1}{\sqrt{2N}}[1, 0, 0, 0, 0]
   \]
   as a tensor product
  \[
  \Theta_{L, \ell}^{(1,0)}(\tau,w) = \frac{1}{\sqrt{2N}}\left( \theta_{N}(\tau;1) \otimes v^2\overline{\theta_{3}(\tau;0)}^4\right).
  \]
  Then by \eqref{eq:etaid} for each $g \in M^!_{\frac{3}{2}, \Zb[N]}$, we have
  \begin{align*}
    \left( g(\tau), \theta_N(\tau;1) \right)^\reg
    &=
      \frac{1}{16}    \left( g(\tau) \otimes( \eta^{-4}(\tau) \vf_4), 
      \theta_N(\tau;1) \otimes v^2\overline{\theta_{3}(\tau;0)}^4  \right)^\reg = \frac{\sqrt{2N}}{16} \Phi_{L, \ell}^{(1,0)}(f,w),
  \end{align*}
  where $f :=  g \otimes( \eta^{-4} \vf_4) \in M_{-\frac{1}{2}, L}$ with $\vf_4$ defined in \eqref{eq:etaid}. From Proposition~\ref{prop:fourier expansion isotropic}, we see that
  \begin{align*}
    \left( g(\tau), \theta_N(\tau,1) \right)^\reg
    &= \frac{\sqrt{2N}}{16} \Phi_{L, \ell}^{(1,0)}(f,w) = \frac{\sqrt{N}}{16}  \{ f, F_{\ell,w}\}
  \end{align*}
  with the power series $F_{\ell,w}(\tau) = \sum_{m, h}  c(h, m) q^m \ef_h \in V_{\frac{5}{2}, L^-}$ defined by
  \begin{equation}
    \label{eq:Flw}
          c(h, m) = -2\sum_{\substack{\lambda  \in K + p(h) \\ -Q(\lambda) = m \geq 0}}
\mathbb{B}_1\left(\frac{p_w(\lambda)}{|\ell_w|} + (h, \ell')\right)+ \sum_{\substack{\lambda \in L + h \\ -Q(\lambda) = m > 0 \\ (\lambda, \ell) \neq 0}} (\sgn(p_w(\lambda)) - \sgn((\lambda, \ell))),
  \end{equation}
  where we understand that the first sum vanishes if $h \notin L_0'/L$.
By \eqref{eq:Bsim}, we have $\frac{p_w(\lambda)}{|\ell_w|} + (h, \ell') \in  \frac{N_3}{6N}$. 
This implies that $F_{\ell,w}$ has coefficients in $\frac{N_3}{3N}\Zb$.
  
If we write $F_{\ell, w} = \sum_{\alpha \in A_N,~ \mu \in A_3^4} F_{\ell, w, \alpha, \mu}$, then the Laurent series
  \begin{equation}
    \label{eq:ttheta}
    \phi^+_N(\tau) := \frac{1}{16} \eta^{-4}(\tau)
    \sum_{\alpha \in A_N} \ef_\alpha
    \sum_{\mu \in A_3^4} \langle\ef_\mu, \vf_4 \rangle
    F_{\ell,w, \alpha, \mu}(\tau)
\in V_{\frac{1}{2},\Zb[N]^-}    
  \end{equation}
  has coefficients in $\frac{N_3}{48N}\Zb$ and satisfies 
  \[
  \left\{ g(\tau), \phi^+_N(\tau) \right\} = \left( g(\tau), \frac{1}{\sqrt{N}} \theta_N(\tau;1)\right)^\reg
  \]
  for any $g \in M_{\frac{3}{2},\Zb[N]}^!$. Hence, Proposition~\ref{prop:xipreimage} implies that $\phi^+_N(\tau)$ is the holomorphic part of $\phi \in S(N)$.  Moreover, $48(N/N_3)\phi^+$ has integral coefficients, which means $  \db(N) \mid 48 N /N_3    $.

  We can now repeat the argument using the lattices $L_j := \Zb[N] \oplus P_j^-$ for $j = 1, 2$ with $P_1, P_2$ the unique two Niemeier lattices with 48 and 72 norm 2 vectors respectively.
Since $\theta_{P_2}(\tau) - \theta_{P_1}(\tau) = 24 \Delta(\tau)$,  we can express
  $$
  \left( g(\tau), \theta_N(\tau;1) \right)^\reg
 = \frac{\sqrt{2N}}{24} \lp \Phi_{L_2, \ell}^{(1,0)}(f,w) - \Phi_{L_1, \ell}^{(1,0)}(f,w)\rp
 $$
 with $f(\tau) = g(\tau) / \Delta(\tau)$. The same argument then shows that $\db(N) \mid 24N$ for all $N \in \Nb$. Combining together with the bound $\db(N) \mid 48N/N_3$ gives us
  \begin{equation}
    \label{eq:bound1}
  \db(N) \mid 24 N /N_3 
  \end{equation}
  for all $N \in \Nb$.

  To obtain a bound without dependence on $N$, we apply the bound \eqref{eq:DMZbound} from \cite{DMZ12} and the Hecke operators in section \ref{subsec:Hecke-op}.
   For any $N \in \Nb$, we write $N = 2^a3^b N'$ with $N'_6 = 1$.
If $a, b \le 1 $, then \eqref{eq:DMZbound} and \eqref{eq:bound1} implies that $\db(N) \mid 24N_2$.
In general, Lemma \ref{lemma:Hecke-op} implies that
\begin{equation}
  \label{eq:bound2}
  \db(Np^e) \mid 2 \lcm(\db(N), \db(Np))
\end{equation}
for all $N \in \Nb, e \ge 1$ and prime $p \nmid N$. 
Therefore when $a \le 1$ and $b \ge 2$, we have
$$
\db(N) \mid \gcd(24N_2\cdot 3^{b}, 2 \lcm(\db(N_2 N'), \db(3N_2 N')))
\mid \gcd(24N_2\cdot 3^b, 48N_2)= 24N_2.
$$
When $a \ge 2$, we have
$$
\db(N) \mid 2 \lcm(\db(3^bN'), \db(2\cdot 3^b N')) \mid 48N_2 = 24 N_4.
$$
\begin{rmk}
  \label{rmk:3N}
  If $3N$ is the sum of two squares, then we can carry out the argument above with $L = \Zb[N] \oplus (\Zb[3]^-)^2$ to obtain $\db(N) \mid 12N_4$.
\end{rmk}

\section{Proof of Theorem~\ref{thm:wt1} and Example \ref{exmp:eta2}}
\label{sec:proof theorem weight one}
For any integral $\mathcal{O}_D$-ideal $\af$ co-prime to $DM$, denote $R = L_{\af, M}$ an anisotropic lattice of signature $(1,1)$ as in defined in \eqref{eq:LaM}.
Since the eigenform $f_\varphi$ is a linear combination of components of the vector-valued cusp form
\begin{equation}
  \label{eq:vartheta}
  \vartheta_R (\tau) := \sum_{\substack{\lambda \in \Gamma_R \backslash R' \\ Q(\lambda) > 0}}
   \sgn(\lambda) q^{Q(\lambda)} \ef_\lambda \in S_{1, R},
\end{equation}
it suffices to consider its $\xi$-preimage.
This is 
constructed in (5.2.1) in \cite{CL20} as
\[
\widetilde{\vartheta}_R(\tau) :=
\log \varepsilon_R \cdot \widetilde{\Theta}_R(\tau) + 
2I'(\tau, R^-),
\]
where $I'$ is a deformed theta integral and $\widetilde{\Theta}_R(\tau)$ is a real-analytic modular form of weight one satisfying
  \begin{equation}
    \label{eq:xi1}
   \xi_1 \widetilde{\Theta}_R(\tau) = \sqrt{2} \Theta_R^{(1,0)}(\tau, z_0),
 \end{equation}
 where  $z_0 = \frac{1}{\sqrt{2}}[1,1]$.
 Here $\Theta_R^{(1,0)}(\tau, z)$ denotes the theta function~\eqref{eq:theta function}. Using the lattice embedding $\iota$ in \eqref{eq:iota} and the map in \eqref{eq:sublat} induced by $\iota$, which we also denote by $\iota$, we can write
 \begin{equation}
   \label{eq:ThetaR}
\sqrt{2}   \Theta_R^{(1,0)}(\tau, z_0)
   = \frac{1}{2\sqrt{AMD}} \lp v^{1/2} \theta_{AMD}(\tau; 1) \otimes\overline{\theta_{AM}(\tau; 0)} \rp \circ \iota.
 \end{equation}
 
As in the proof of Theorem 5.1 loc.\ cit., the number $\kappa_R$ defined in (4.2.7) is a bound on the denominator  of Fourier coefficients of the holomorphic part\footnote{The holomorphic function $\Theta^+_R(\tau)$ is also called a mixed mock modular form in the sense of \cite{DMZ12}.}
  $\Theta^+_R(\tau)$ of $\widetilde{\Theta}_R(\tau)$.
For this purpose, it is enough to produce a suitable $\Theta^+_R(\tau)$ having rational Fourier coefficients with a good denominator bound. 
To do this, we can simply take $\Theta_R^+(\tau) = \frac{1}{2}(\phi^+(\tau) \otimes \theta_{AM}(\tau; 0))\circ \iota$ for any $\phi \in S(AMD)$ as in the proof of Theorem \ref{thm:denom}. 
Since we can choose $A$ to be odd, this gives the bound $\kappa_R = 48(DM)_4$, which leads to the improvement of $\kappa = 96(DM)_4$ following the same proofs of Theorems 5.1 and 6.5 in \cite{CL20}. This finishes the proof of Theorem~\ref{thm:wt1}.

For Example \ref{exmp:eta2}, we take $F = \Qb(\sqrt{12})$ and $\af = \Oc_F, M = 2$. Then $R' = (2\sqrt{3})^{-1}\Zb[\sqrt{3}]$, $A_R = \{\frac{a + b\sqrt{3}}{2\sqrt{3}}: a \in \Zb/12\Zb, b \in \Zb/4\Zb\} \cong A_6 \oplus A_2$ and $\varepsilon_R = (2 + \sqrt{3})^4$.
Furthermore, we have
\begin{equation}
  \label{eq:eta2}
  \begin{split}
    \langle\vartheta_R(\tau), \vf\rangle &= \eta(\tau)^2 , \qquad
\langle \Theta_R^{(1, 0)}(\tau, z_0), \vf \rangle = \frac{1}{4} \eta(\tau)^3 \eta(-\overline{\tau}) \sqrt{v},\\
    \vf &:= \frac{1}{8}\vf_6 \otimes \vf_4 \in \Qb[A_6] \otimes \Qb[A_2] \cong \Qb[A_R].    
  \end{split}
\end{equation}
Let $\widetilde{\Theta}^+_R$ be as above and $\widetilde{\Theta}_R$ its modular completion.
Then $\langle \widetilde{\Theta}^+_R, \vf \rangle$ is in $q^{-1/12}\Qb\llbracket q\rrbracket$ and the $\xi_1$-image of its modular completion  $\langle \widetilde{\Theta}_R, \vf \rangle$ is $\eta(\tau)^3 \eta(-\overline{\tau})\sqrt{v}/(2\sqrt{2})$.
From \eqref{eq:mocketa3} and Remark \ref{rmk:chir}, we then know that $\eta(\tau) \widetilde{\theta}^+_2 - \langle \widetilde{\Theta}^+_R, \vf \rangle$ is in the trivial space $M_{1, \overline{\chi}^2}$ with $\chi$ the character defined in \eqref{eq:chi}. 
Finally, Proposition 5.5 in \cite{CL20} shows that the holomorphic part of 
$\langle 2 I'(\tau, R^-), \vf \rangle$
is given by
$\sum_{\af \subset \Oc_F} \widetilde{\varphi}(\af) q^{\Nm(\af)/12} \in q^{11/12} \Rb\llbracket q \rrbracket$.
Therefore, $\widetilde{\vartheta}^+$ defined in \eqref{eq:tvartheta} is the holomorphic part of the harmonic Maass form 
$ \langle 2I'(\tau, R^-) + \log \varepsilon_R \widetilde{\Theta}_R(\tau) , \vf \rangle \in H_{1, \overline{\chi}^2}$, whose $\xi_1$-image is $\langle \vartheta_R(\tau), \vf \rangle = \eta(\tau)^2 \in S_{1, \chi^2}$. 

\section{Explicit constructions of mock modular forms}
\label{sec:explicit construction}

In this section we compute explicit mock modular forms of weight $\frac{1}{2}$ and $\frac{3}{2}$ using the evaluations of theta lifts for lattices of signature $(1,1)$ given in Section~\ref{sec:theta lifts}.
The formulas involve the holomorphic (quasimodular) Eisenstein series $E_2(\tau) = 1-24\sum_{n \geq 1}\sigma_1(n)q^n$, and the periodic Bernoulli polynomials $\mathbb{B}_k(x)$ defined in Section~\ref{theta lifts}.

\subsection{Mock modular forms of weight $\frac{1}{2}$}

We construct mock modular forms $\widetilde{\theta}_N^+(\tau;1)$ of weight $\frac{1}{2}$ for $\overline{\rho}_N$ with shadow $\frac{1}{\sqrt{N}}\theta_N(\tau;1)$ for every $N \in\Nb$.

\begin{prop}\label{prop:mockthetaweightonehalf}
	 \begin{enumerate}
		\item Suppose that $2N$ is a square. Then
		\begin{align*}
		\widetilde{\theta}_N^+(\tau;1) &= \frac{1}{\eta(\tau)^3}\bigg(\sum_{\substack{x,y \in \Zb \\ y > \frac{\sqrt{2}|x|}{\sqrt{N}}}}\left( \frac{-4}{y}\right)\sgn(x)yq^{\frac{y^2}{8}-\frac{x^2}{4N}} \ef_{x} \\
		&\qquad \qquad \qquad + \sum_{b(2\sqrt{2N})}\left(\frac{-4}{b}\right)\left(\frac{E_2(\tau)}{12\sqrt{2N}}-\sqrt{2N}\mathbb{B}_2\left(\frac{b}{2\sqrt{2N}} \right) \right)\ef_{b\sqrt{\frac{N}{2}}}\bigg)
		\end{align*}
		is a mock modular form of weight $\frac{1}{2}$ for $\overline{\rho}_N$ with shadow $\frac{1}{\sqrt{N}}\theta_{N}(\tau;1)$. 
		\item Suppose that $2N$ is not a square. Let $F = \Qb(\sqrt{2N})$ and let $\varepsilon_N = a+b\sqrt{2N}$ be the smallest totally positive unit $> 1$ of $F$ such that $b$ is even and $\lcm(2N, 4) \mid (a-1)$. Then 
		\[
		\widetilde{\theta}_N^+(\tau;1) = \frac{1}{\eta(\tau)^3}\sum_{\substack{x,y \in \Zb  \\ \varepsilon_N^{-2} < \frac{\sqrt{N}y+\sqrt{2}x}{\sqrt{N}y-\sqrt{2}x}< 1}}\left( \frac{-4}{y}\right)\sgn(x)
		\Tr_{F/\Qb}\left(\frac{\frac{y}{2}+\frac{x}{\sqrt{2N}}}{1-\varepsilon_N^{-1}}\right)
		q^{\frac{y^2}{8}-\frac{x^2}{4N}} \ef_{x}
		\]
		is a mock modular form of weight $\frac{1}{2}$ for $\overline{\rho}_N$ with shadow $\frac{1}{\sqrt{N}}\theta_{N}(\tau;1)$. 
	\end{enumerate}
\end{prop}

\begin{proof}
	The proof is similar to the proof of Theorem~\ref{thm:denom}, but we will give some details for the convenience of the reader. Let us first assume that $2N$ is a square. Then the lattice
	\[
	L = \sqrt{N}[1,1]\Zb \oplus \sqrt{2}[1,-1]\Zb
	\]
	in $(\Rb^2,xy)$ has signature $(1,1)$ and is isotropic. Its dual lattice is given by
	\[
	L' = \frac{1}{2\sqrt{N}}[1,1]\Zb \oplus \frac{1}{2\sqrt{2}}[1,-1]\Zb.
	\]
	We choose the primitive isotropic vector $\ell = [2\sqrt{N},0] \in L$. Then we have $(\ell,L) = 2\sqrt{2N}$. The theta function $\Theta_{L, \ell}^{(1,1)}(\tau,z)$ considered in Section~\ref{sec:theta lifts} splits at the special point $w = \frac{1}{\sqrt{2}}[1,1]$ as a tensor product
	\[
	\Theta_{L, \ell}^{(1,1)}(\tau,w) = -\frac{1}{2\sqrt{2N}}\theta_{N}(\tau;1) \otimes v^{\frac{3}{2}}\overline{\theta_2(\tau;1)}.
	\]
	Hence, using Lemma~\ref{lemma:eigenvector} we can write for any $g \in M_{\frac{3}{2},\Zb[N]}^!$
	\begin{align*}
	(g(\tau),\theta_N(\tau;1))^\reg &= \frac{1}{2}\left(g(\tau) \otimes \eta^{-3}(\tau) \vf_2,\theta_{N}(\tau;1) \otimes v^{\frac{3}{2}}\overline{\theta_2(\tau;1)}\right)^\reg = -\sqrt{2N} \Phi_{L, \ell}^{(1,1)}(f,w),
	\end{align*}
	where $f = g \otimes (\eta^{-3}\vf_2) \in M_{0,L}^!$. On the other hand, Proposition~\ref{prop:fourier expansion isotropic} implies that the theta lift can be expressed in terms of the pairing $\{\cdot,\cdot\}$ defined in \eqref{eq:pairing} as
	\[
	\Phi_{L, \ell}^{(1,1)}(f,w) = \{f,F_{\ell,w}\}
	\]
	with the power series
	\begin{align*}
	F_{\ell,w}(\tau) &= -\sum_{b(2\sqrt{2N})}\left(\frac{\sqrt{2}E_2(\tau)}{24|\ell_{w^\perp}|}-\frac{|\ell_w|^2}{\sqrt{2}|\ell_{w^\perp}|}\mathbb{B}_2\left(\frac{b}{2\sqrt{2N}}\right)\right)\ef_{b\ell/2\sqrt{2N}} \\
	&\quad + \frac{1}{\sqrt{2}}\sum_{\substack{\lambda \in L' \\ Q(\lambda) < 0}}p_{w^\perp}(\lambda)\big(\sgn(p_w(\lambda))-\sgn((\lambda,\ell))\big)q^{-Q(\lambda)}\ef_{\lambda+L}.
	\end{align*}
	Here we used that the lattice $K$ is trivial, and that a system of representatives for $\delta \in L_0'/L$ is given by $b \ell/M$, where $b$ runs modulo $M:= (\ell,L)$. In particular, it follows from Proposition~\ref{prop:xipreimage} that the Laurent series
	\[
	\widetilde{\theta}_N^+(\tau;1) := -\sqrt{2}\eta^{-3}(\tau)\langle F_{\ell,w}(\tau),\vf_2\rangle \in V_{\frac{1}{2},\Zb[N]^-}
	\]
	is a mock modular form of weight $\frac{3}{2}$ with shadow $\frac{1}{\sqrt{N}}\theta_{N}(\tau;1)$. Using $|\ell_w| = |\ell_{w^\perp}| = \sqrt{2N}, p_w(\lambda) = \frac{x}{\sqrt{2N}}$ and $p_{w^\perp}(\lambda) = -\frac{y}{2}$ for $\lambda = [\frac{x}{2\sqrt{N}}+\frac{y}{2\sqrt{2}},\frac{x}{2\sqrt{N}}-\frac{y}{2\sqrt{2}}] \in L'$ it is easy to obtain the representation of $\widetilde{\theta}_N^+(\tau;1)$ given in the proposition. We leave the details of the simplification to the reader.
	
	The proof for $2N$ not being a square is similar, so we will only give a sketch. In this case, we use a certain lattice $L_{\mathfrak{a},M}$ of the shape \eqref{eq:LaM}. We let $D = 2N$ if $N$ is even and $D = 8N$ if $N$ is odd. We choose
\[
\mathfrak{a} := \begin{cases}
\mathcal{O}_D & \text{if $N$ is even}, \\
2\Zb + \sqrt{2N}\Zb & \text{if $N$ is odd},
\end{cases} \qquad
M := \begin{cases}
2 & \text{if $N$ is even}, \\
1 & \text{if $N$ is odd}.
\end{cases}
\]
Note that in the first case we have $A = [\mathcal{O}_D:\mathfrak{a}] = 1$ and in the second case we have $A= 2$. In any case, we have
\[
(L_{\mathfrak{a},M},Q_{\mathfrak{a},M}) =\left(2\Zb + \sqrt{2N}\Zb ,-\frac{\Nm_{F/\Qb}}{2} \right)  .
\]
We can now proceed similarly as in the case of $2N$ being a square, now using the evaluation of the theta lift $\Phi_{L, \ell}^{(1,1)}(f,w)$ at $w = \frac{1}{\sqrt{2}}[1,1]$ given in Corollary~\ref{cor:unfolding anisotropic}. It is easy to see that the unit $\varepsilon_N:=\varepsilon_{L_{\mathfrak{a},M}}$ is characterized by the conditions given in the proposition. 
\end{proof}

\begin{rmk}
  \label{rmk:1}
	\begin{enumerate}
		\item For $2N$ being a square, in \cite{CL20} the authors constructed the mock modular form
		\begin{align*}
		& \frac{1}{\eta(\tau)^3}\bigg(\sum_{\substack{x,y \in \Zb \\ y > \frac{\sqrt{2}|x|}{\sqrt{N}}}}\left( \frac{-4}{y}\right)\sgn(x)\left(y-\frac{\sqrt{2}|x|}{\sqrt{N}}\right)q^{\frac{y^2}{8}-\frac{x^2}{4N}} \ef_{x} + \frac{E_2(\tau)}{6\sqrt{2N}}\sum_{b(2\sqrt{2N})}\left(\frac{-4}{b}\right)\ef_{b\sqrt{\frac{N}{2}}}\bigg)
		\end{align*}
		of weight $\frac{1}{2}$ with shadow $\frac{1}{\sqrt{N}}\theta_{N}(\tau;1)$. The difference with the mock modular form constructed in Proposition~\ref{prop:mockthetaweightonehalf} is given by $\eta^{-3}$ times an Eisenstein series of weight $2$ for $\overline{\rho}_N$, which can for example be constructed as an integral of a Kudla-Millson theta function as in Theorem 6.4 in \cite{funkemillson11}.
      \item For fixed $m \in \Nb$ the number of $x,y \in \Zb$ with $0 \le Ny^2-2x^2 \le m$ and 
	\[
	\varepsilon_N^{-2} < \frac{\sqrt{N}y+\sqrt{2}x}{\sqrt{N}y-\sqrt{2}x} < 1
	\]
	is finite. Indeed, the above conditions imply that $x,y$ satisfy the inequalities
	\[
	\frac{\sqrt{m}}{2\sqrt{N}}(1+\varepsilon_N^{-1})<|y| < \frac{\sqrt{m}}{2\sqrt{N}}(1+\varepsilon_N), \qquad \frac{\sqrt{m}}{2\sqrt{2}}(\varepsilon_N^{-1}-\varepsilon_N) < \sgn(y)x < 0.
	\]
	In particular, the coefficient at $q^m$ in the inner sum in $\widetilde{\theta}_N^+(\tau;1)$ in item (2) of Proposition~\ref{prop:mockthetaweightonehalf} is given by a finite sum of rational numbers. 
	\item If $2N$ is a square, then the denominators of the coefficients of $\widetilde{\theta}_N^+(\tau;1)$ are bounded by $6\sqrt{2N}$. If $2N$ is not a square, then we can write 
	\[
	\Tr_{F/\Qb}\left(\frac{\lambda}{1-\varepsilon_N^{-1}}\right) = \frac{\Tr_{F/\Qb}(\lambda(1-\varepsilon_N))}{\Nm_{F/\Qb}(1-\varepsilon_N^{-1})} = \frac{\Tr_{F/\Qb}(\lambda(1-\varepsilon_N))}{\Tr_{F/\Qb}(1-\varepsilon_N)},
	\]
	which implies that the denominators of the coefficients of $\widetilde{\theta}_N(\tau;1)$ are bounded by $\Tr_{F/\Qb}(1-\varepsilon_N)$. In the numerical examples that we looked at the denominators of the coefficients of $\widetilde{\theta}_N(\tau;1)$ were usually bigger than the bound from Theorem~\ref{thm:denom}. 
	\end{enumerate}
\end{rmk}

\begin{exmp}\label{example mock theta order 3}
	Consider Ramanujan's classical mock theta functions of order $3$,
	\begin{align*}
	f(q) &= 1+ \sum_{n=1}^\infty \frac{q^{n^2}}{(1+q)^2(1+q^2)^2 \cdots (1+q^n)^2} =  1+q-2q^2 + 3q^3 - 3q^4 + 3q^5 + \dots\\
	\omega(q) &= \sum_{n=0}^\infty \frac{q^{2n(n+1)}}{(1-q)^2(1-q^3)^2 \cdots(1-q^{2n+1})^2} = 1 + 2q + 3q^2 + 4q^3 + 6q^4 + 8q^5 + \dots
	\end{align*}
	The fundamental work of Zwegers \cite{ZwegersPaper} (see also Section~8.2 in \cite{bruinierono}) shows that the $\Cb[\Zb/12\Zb]$-valued function 
	\begin{align*}
	F^+(\tau) &= q^{-\frac{1}{24}}f(q) (\ef_1 - \ef_5 + \ef_7 - \ef_{11})\\
			& \quad + 2q^{\frac{1}{3}}(\omega(q^{\frac{1}{2}})-\omega(-q^{\frac{1}{2}}))(-\ef_2 + \ef_{10}) \\
			& \quad + 2q^{\frac{1}{3}}(\omega(q^{\frac{1}{2}})+\omega(-q^{\frac{1}{2}}))(-\ef_4 + \ef_{8})
	\end{align*}
	is the holomorphic part of a harmonic Maass form $F(\tau)$ of weight $\frac{1}{2}$ for the dual Weil representation associated with $\Zb[6]$. By comparing principal parts, we obtain that
	\[
	F^+(\tau) = 2\widetilde{\theta}_6(\tau;1)+2\widetilde{\theta}_6(\tau;1)^{\sigma_3},
	\]
	where $\widetilde{\theta}_6(\tau;1)$ is the mock modular form constructed in the second item of Proposition~\ref{prop:mockthetaweightonehalf} above, and $\sigma_3$ is the Atkin-Lehner involution on $\Zb/12\Zb$ that interchanges $1$ with $7$, $3$ with $9$, $5$ with $11$, and fixes all other elements. For example, this implies the identity
	\begin{align*}
	q^{-\frac{1}{24}}f(q) &= 2 \widetilde{\theta}_{6,1}(\tau;1) + 2 \widetilde{\theta}_{6,7}(\tau;1)\\
	& = \frac{1}{\eta(\tau)^3}\sum_{\substack{x,y \in \Zb \\ x \equiv 1 (6) \\ \varepsilon_6^{-2} < \frac{\sqrt{3}y+x}{\sqrt{3}y-x}< 1}}\left( \frac{-4}{y}\right)\sgn(x)
		\Tr_{F/\Qb}\left(\frac{y+\frac{x}{\sqrt{3}}}{1-\varepsilon_6^{-1}}\right)
		q^{\frac{y^2}{8}-\frac{x^2}{24}},
	\end{align*}
	where $F = \Qb(\sqrt{12})$ and $\varepsilon_6 = 97+28\sqrt{12}$.
\end{exmp}

If $6N$ is a square then we can use the theta lift $\Phi^{(1, 0)}_{L, \ell}(f,z)$ 
on the isotropic lattice $L = \sqrt{N}[1,1]\Zb \oplus \sqrt{6}[1,-1]\Zb$, and apply similar arguments as in the proof of Proposition~\ref{prop:mockthetaweightonehalf} to obtain the following formula.

\begin{prop}\label{prop:mock theta weight onehalf alternative}
 Suppose that $6N$ is a square. Then
	\[
	\widetilde{\theta}_N^+(\tau;1) = \frac{1}{\eta(\tau)}\bigg(\sum_{\substack{x,y \in \Zb \\ y >\frac{\sqrt{6}|x|}{\sqrt{N}}}}\left(\frac{12}{y}\right)\sgn(x)q^{\frac{y^2}{24}-\frac{x^2}{4N}}\ef_x-\sum_{b(2\sqrt{6N})}\left( \frac{12}{b}\right)\mathbb{B}_1\left(\frac{b}{2\sqrt{6N}}\right)\ef_{b\sqrt{\frac{N}{6}}} \bigg)
	\]
	is a mock modular form of weight $\frac{1}{2}$ for $\overline{\rho}_N$ with shadow $\frac{1}{\sqrt{N}}\theta_{N}(\tau;1)$.
\end{prop}

\begin{exmp}
	 Arguing as in Example~\ref{example mock theta order 3}, we obtain from Proposition~\ref{prop:mock theta weight onehalf alternative} the identity
	\[
	q^{-\frac{1}{24}}f(q) = \frac{1}{\eta(\tau)}\bigg(1 + 2\sum_{\substack{x,y \in \Zb \\ x \equiv 1 (6) \\ y > |x|}}\left(\frac{12}{y}\right)\sgn(x)q^{\frac{y^2}{24}-\frac{x^2}{24}} \bigg)
	\] 
	for Ramanujan's order $3$ mock theta function $f(q)$.
	This identity is very similar to the ones for Ramanujan's mock theta functions of order $5$ and $7$ given in Section 6 of \cite{Za09}.
\end{exmp}

\begin{rmk}
	Unfortunately, if $6N$ is not a square, our method does not work as before. In this case we would need to evaluate the theta lift $\Phi_{L, \ell}^{(1,0)}(f,z)$ on an anisotropic lattice $L$ of signature $(1,1)$, where $f$ is a weight $1$ weakly holomorphic modular form for $\rho_L$. Since there might be holomorphic modular forms of weight $1$ for $\rho_L$, we cannot write $f$ as a linear combination of Maass Poincar\'e series to compute the theta lift of $f$ by the unfolding argument as in Proposition~\ref{prop:unfolding anisotropic}. We wonder if our method can be adjusted to obtain a mock modular form with shadow $\frac{1}{\sqrt{N}}\theta_N(\tau;1)$ using the theta lift $\Phi_{L, \ell}^{(1,0)}(f,z)$ for anisotropic $L$.
\end{rmk}

\subsection{Mock modular forms of weight $\frac{3}{2}$}

We now give explicit mock modular forms of weight $\frac{3}{2}$ with shadow $\frac{\sqrt{N}}{\pi}\theta_N(\tau;0)$. The construction works analogously as in the proof of Proposition~\ref{prop:mockthetaweightonehalf}, so we omit the details for brevity. 

Using the signature $(1,1)$ theta lift $\Phi_{L, \ell}^{(0,0)}(f,z)$ with a constant polynomial considered in Section~\ref{sec:theta lifts} we obtain the following mock modular forms.

\begin{prop}\label{prop:mock theta weight 3half}
\begin{enumerate}
	\item Suppose that $6N$ is a square. Then
	\begin{align*}
	\widetilde{\theta}_N^+(\tau;0) &= \frac{1}{\eta(\tau)}\bigg(-2\sum_{\substack{x,y \in \Zb \\ y >\frac{\sqrt{6}|x|}{\sqrt{N}}}}\left(\frac{12}{y}\right)|x|q^{\frac{y^2}{24}-\frac{x^2}{4N}}\ef_x \\ 
	&\qquad\qquad \qquad +\sum_{b(2\sqrt{6N})}\left( \frac{12}{b}\right)\left(\frac{ E_2(\tau)}{12} +2N\mathbb{B}_2\left( \frac{b}{2\sqrt{6N}}\right)\right)\ef_{b\sqrt{\frac{N}{6}}} \bigg)
	\end{align*}
	is a mock modular form of weight $\frac{3}{2}$ for $\overline{\rho}_N$ with shadow $\frac{\sqrt{N}}{\pi}\theta_{N}(\tau;0)$.
	\item Suppose that $6N$ is not a square. Let $F = \Qb(\sqrt{6N})$ and let $\varepsilon_N = a+b\sqrt{6N}$ be the smallest totally positive unit $> 1$ of $F$ such that $b$ is even and $\lcm(2N, 12) \mid (a-1)$. Then
	\[
	\widetilde{\theta}_N^+(\tau;0) = \frac{1}{\eta(\tau)}\sum_{\substack{x,y \in \Zb \\ \varepsilon_N^{-2} < \frac{\sqrt{N}y+\sqrt{6}x}{\sqrt{N}y-\sqrt{6}x} \leq 1}}\left( \frac{12}{y}\right)\sgn(y)
		\Tr_{F/\Qb}\left(\frac{\frac{\sqrt{N}}{\sqrt{6}}y+x}{1-\varepsilon_N^{-1}}\right)
		q^{\frac{y^2}{24}-\frac{x^2}{4N}} \ef_{x}
		\]
		is a mock modular form of weight $\frac{3}{2}$ for $\overline{\rho}_N$ with shadow $\frac{\sqrt{N}}{\pi}\theta_{N}(\tau;0)$.
\end{enumerate}
\end{prop}

Similarly, using the signature $(1,1)$ theta lift $\Phi_{L, \ell}^{(0,1)}(f;z)$ with a degree $(0,1)$ polynomial, we obtain the following result.

\begin{prop} 
\label{prop:mock theta weight 3half alternative}
\begin{enumerate}
	\item Suppose that $2N$ is a square. Then
	\begin{align*}
	\widetilde{\theta}_N^+(\tau;0) &= \frac{1}{\eta^3(\tau)}\bigg(-2\sum_{\substack{x,y \in \Zb \\ y >\frac{\sqrt{2}|x|}{\sqrt{N}}}}\left(\frac{-4}{y}\right)|x|yq^{\frac{y^2}{8}-\frac{x^2}{4N}}\ef_x \\ 
	&\qquad\qquad \qquad +\frac{16N^2}{3\sqrt{2N}}\sum_{b(2\sqrt{2N})}\left( \frac{-4}{b}\right)\mathbb{B}_3\left( \frac{b}{2\sqrt{2N}}\right)\ef_{b\sqrt{\frac{N}{2}}} \bigg)
	\end{align*}
	is a mock modular form of weight $\frac{3}{2}$ for $\overline{\rho}_N$ with shadow $\frac{\sqrt{N}}{\pi}\theta_{N}(\tau;0)$. 
	\item Suppose that $2N$ is not a square. Let $F = \Qb(\sqrt{2N})$ and let $\varepsilon_N = a+b\sqrt{2N}$ be the smallest totally positive unit $> 1$ of $F$ such that $b$ is even and $\lcm(2N, 4) \mid (a-1)$. Then
	\[
	\widetilde{\theta}_N^+(\tau;0) = \frac{1}{\eta(\tau)^3}\sum_{\substack{x,y \in \Zb \\ \varepsilon_N^{-2} < \frac{\sqrt{N}y+\sqrt{2}x}{\sqrt{N}y-\sqrt{2}x} \leq 1}}\left( \frac{-4}{y}\right)\sgn(y)
		\Tr_{F/\Qb}\left(\frac{\sqrt{2N}(\frac{y}{2}+\frac{x}{\sqrt{2N}})^2}{1-\varepsilon_N^{-2}}\right)
		q^{\frac{y^2}{8}-\frac{x^2}{4N}} \ef_{x}
		\]
		is a mock modular form of weight $\frac{3}{2}$ for $\overline{\rho}_N$ with shadow $\frac{\sqrt{N}}{\pi}\theta_{N}(\tau;0)$. 
\end{enumerate}
\end{prop}

\bibliographystyle{plain}

\bibliography{unary}{}

\end{document}